\documentclass[a4paper,twoside]{article}
\usepackage{amsmath,amsfonts,amssymb,amstext,amssymb}
\usepackage[english]{babel}
\usepackage[utf8]{inputenc}
\usepackage{hyperref}
\hypersetup{
    bookmarks=true,         
    unicode=true,          
    pdftoolbar=true,        
    pdfmenubar=true,        
    pdffitwindow=true,      
    pdftitle={},    
    pdfauthor={},     
    pdfsubject={},   
    pdfcreator={},   
    pdfproducer={}, 
    pdfkeywords={}, 
    colorlinks=true,       
    linkcolor=blue,  
    citecolor=blue,        
    filecolor=magenta,      
    urlcolor=violet           
}
\usepackage{array}
\usepackage{theorem}
\usepackage{enumerate}
\usepackage{stmaryrd}
\usepackage{float}
\usepackage{graphics,epsfig,times,epstopdf}
\usepackage{ae,aecompl}
\usepackage{mathrsfs}
\usepackage[all]{xy}
\usepackage{graphicx}
\usepackage{epstopdf}
\usepackage{fancyhdr}
\usepackage{here}
\usepackage{dsfont}
\usepackage{caption}
\usepackage{tikz}
\usetikzlibrary{positioning,decorations,arrows,patterns,shapes,backgrounds}

\theoremstyle{plain}
\theoremstyle{definition}

\newtheorem{remark}[equation]{Remark}

\newtheorem{definition}[equation]{Definition}
\numberwithin{equation}{section}
\newcommand{\OO}{\mathcal {O}}

\renewcommand{\P}{\mathcal{P}}
\newcommand{\N}{\mathcal{N}}

\newcommand{\bo}{\rm bo}
\newcommand{\Bo}{\rm Bo}

\newcommand{\mfT}{\mathfrak T}
\newcommand{\dsp}{\displaystyle}
\newcommand{\RR}{\mathbb{R}}
\newcommand{\NN}{\mathbb{N}}

\newcommand{\zetabar}{\underline{\zeta}}
\newcommand{\psibar}{\underline{\psi}}
\newcommand{\vbar}{\underline{v}}
\newcommand{\id}[1]{\left\vert_{\scriptstyle #1}\right.}
\newcommand{\nn}{\nonumber}
\renewcommand{\t}{\tilde}
\usepackage{authblk}
\title{A numerical scheme for an improved Green–Naghdi model in the Camassa-Holm regime for the propagation of internal waves}
\author[1]{{Christian Bourdarias} \thanks{\url{christian.bourdarias@univ-smb.fr}}}
\author[1]{{St\'ephane Gerbi}\thanks{\url{stephane.gerbi@univ-smb.fr}}}
\author[1,2]{{Ralph Lteif}\thanks{\url{ralphlteif_90@hotmail.com}}}
\affil[1]{\small{Laboratoire de Math\'ematiques, UMR 5127 - CNRS and Universit\'e Savoie Mont Blanc\par73376 Le Bourget-du-Lac Cedex, France.}}
\affil[2]{\small{Laboratory of Mathematics-EDST, Lebanese University, Beirut, Lebanon. }}
\date{}
\begin{document}
\maketitle
 
\begin{abstract}
In this paper we introduce a new reformulation of the Green-Naghdi model in the Camassa-Holm regime for the propagation of internal waves over a flat topography derived by 
Duch\^ene, Israwi and Talhouk [{\em SIAM J. Math. Anal.}, 47(1), 240–-290]. These new Green–Naghdi systems are adapted to improve the frequency dispersion of the original model, they share the same order of precision as the standard one but have an appropriate structure which makes them much more suitable for the numerical resolution. We develop a second order splitting scheme where the hyperbolic part of the system is treated with a high-order finite volume scheme and the dispersive part is treated with a finite difference approach. Numerical simulations are then performed to validate the model.
\end{abstract}
\textbf{Key words} : Green–Naghdi model, Nonlinear shallow water, Dispersive waves, Nonlinear interactions, Improved dispersion, Splitting method,  Hybrid method, Finite volume, Finite Difference, High order scheme, WENO reconstruction

\newpage
\tableofcontents
\newpage
\section{Introduction}
This study deals with the propagation of internal waves in the uni-dimensional setting located at the interface between two layers of fluids of different densities. The fluids are assumed to be incompressible, homogeneous, 
and immiscible, limited from above by a rigid lid and from below by a flat bottom. 
This type of fluid dynamics problem is encountered by researchers in oceanography when they study the wave near the shore. 
Because of the difference in the salinity of the different layers of water near the shore, it is useful to model the flow of salted water by a two layers incompressible fluids flow.
The usual way of describing such a flow is to use the 3D-Euler equations for the different layers adding some thermodynamic and 
dynamic conditions at the interface. This system will be called the \emph {full Euler system}. This system of partial differential
equations is very rich but very difficult to manipulate both mathematically and numerically. 
This is the reason why reduced models have been derived to characterize the evolution of the solution in some physical and geometrical specific regimes. Many models for the water wave (air-water interface) system have already been derived and studied in the shallow-water regime, where we consider 
that the wave length of the flow is large compared to the typical depth. We refer the reader to the following papers \cite{Barthelemy04,BD08,LB09,BBCCLMT11,Lannes}. 
Earlier works have also set a very interesting theoretical background for the two-fluid system see~\cite{BonaLannesSaut08,Anh09,Duchene13,DucheneIsrawiTalhouk14,DucheneIsrawiTalhouk15},  
for more details.

One of these reduced models is the Green-Naghdi system of partial differential equations (denoted GN in the following). Many numerical studies of the  fully nonlinear GN equations for the one layer case have been proposed in the literature. Let us introduce some of these results. The GN model describing dispersive shallow water waves have been numerically studied in~\cite{MGH10} after being written in terms of potential variables in a pseudo-conservative form and using a Godunov scheme. Bonneton {\em et al.} studied numerically in~\cite{BCLMT} the fully nonlinear and weakly dispersive GN model for shallow water waves of large amplitude over variable topography. The original model was firstly adjusted in a suitable way for numerical resolution with an improvement of the dispersive properties, then they proposed to use a finite volume-finite difference splitting scheme that decomposes the hyperbolic and dispersive parts of the equations. This similar approach was earlier introduced in~\cite{EIK05} for the solution of the Boussinesq equations. In~\cite{CLM}, a three-parameters GN system is derived, yielding further improvements of the dispersive properties. Let us mention also~\cite{MID14}, where a highly accurate and stable numerical scheme based on the Galerkin / finite-element method was presented for the Serre system. The method is based on smooth periodic splines in space and an explicit fourth-order Runge-Kutta method in time.
Recently, Lannes and Marche introduced in~\cite{LannesMarche14} a new class of two-dimensional GN equations over varying topography. These fully nonlinear and weakly dispersive equations have a mathematical structure more suitable for 2D simulations. Using the same splitting strategy initiated in~\cite{CLM}, they develop a high order, well balanced, and robust numerical code for these new models. Finally, we would like to mention the recent work of Duch\^ene, Israwi and Talhouk~\cite{DucheneIsrawiTalhouk16}, where they derive a new class of GN models with improved frequency dispersion for the propagation of internal waves at the interface between two layers of fluids. They numerically compute their class of GN models using spectral methods~\cite{Trefethen00} for space discretization and the Matlab solver ode45, which is based on the fourth and fifth order Runge-Kutta-Merson method for time evolution. Their numerical simulations show how the different frequency dispersion of the modified GN models may affect the appearance of high frequency Kelvin-Helmholtz instabilities.

In this paper, we present the numerical resolution of the GN model in the Camassa-Holm (medium amplitude) regime obtained and fully justified by Duch\^ene, Israwi and Talhouk in \cite{DucheneIsrawiTalhouk15} using the same strategy initiated in~\cite{BCLMT,LannesMarche14}. Let us recall that this model describes the propagation of one-dimensional internal waves at the interface between two layers of ideal fluids, limited by a rigid lid from above and a flat bottom from below. This model is first recast under a new formulation more suitable for numerical resolution with the same order of precision as the standard one but with improved frequency dispersion. More precisely, following~\cite{BCLMT} we derive a family of GN equations depending on a parameter $\alpha$ to be chosen in order to minimize the phase velocity error between the reduced model and the \emph {full Euler system}. Then we propose a numerical scheme that decomposes the hyperbolic and dispersive parts of the equations. Following the same strategy adopted in~\cite{BCLMT,LannesMarche14}, we use a second order splitting scheme.
The approximation $U^{n+1}=(\zeta^{n+1}, v^{n+1})$ is computed at time $t^{n+1} = t^n + \Delta t$ where $\zeta$ represents the deformation of the interface and $v$ represent the {\em shear mean velocity} defined in Section~\ref{GNCHsec},
in terms of the approximation $U^n$ at time $t^n$ by solving
$$U^{n+1} = S_1(\Delta t/2)S_2(\Delta t)S_1(\Delta t/2)U^n,$$
where $S_1(.)$ is the operator associated to the hyperbolic part and $S_2(.)$ the operator associated to the dispersive part of the GN equations.
For the numerical computation of $S_1(.)$, we use a finite volume method. We begin by the VFRoe method (see~\cite{BGH00,GHN02,GHN03}), that is an approximate Godunov scheme. Unfortunately, the VFRoe scheme seems to be very diffusive. To this end, we propose a second-order scheme following the classical ``MUSCL" approach~\cite{Leer79}. Finally, following~\cite{JiangShu96} we implement a fifth-order accuracy WENO reconstruction in order to reach higher order accuracy in smooth regions and a good resolution around discontinuities since the second order schemes are known to degenerate to first
order accuracy at smooth extrema. On the other hand, $S_2(.)$ is computed using a finite difference scheme discretized using second and fourth order formulas, whereas for time discretization we use classical second and fourth order Runge-Kutta methods according to the order of the space derivative in consideration.  
 The computation of $S_2(\Delta t)$ in the above splitting scheme (dispersive part) requires the inversion of the following symmetric second order differential operator introduced in \cite{DucheneIsrawiTalhouk15}
$$
      {\mathfrak T}[\epsilon\zeta]V \ = \ q_1(\epsilon\zeta)V \ - \ \mu \nu \partial_x \Big(q_2(\epsilon\zeta)\partial_xV \Big).
$$
with $q_i(X)\equiv 1+\kappa_i X $ ($i=1,2$), where $\kappa_i$ and $\nu$ are constants depending on three parameters: the ratio between the densities of the two layers, the ratio between their depth and the capillary effect of the interface.   
Moreover, it is known that third order derivatives involved in this model may create high frequencies instabilities, but the presence of the operator ${\mathfrak T}[\epsilon\zeta]^{-1}$ in the second equation stabilizes the model with respect to these perturbations, allowing for more robust numerical computations, see Section~\ref{IHFsec}. The invertibility of this operator played also an important role in the well-posedness of these equations (see~\cite{Israwi11} for the one layer case and~\cite{DucheneIsrawiTalhouk15} for the two layers case).      
However, this operator is time dependent since it depends on the deformation of the interface $\zeta(t,x)$. In order to avoid the inversion of this operator at each time step keeping its stabilizing effects, we derive a new family of physical models that are equivalent to the standard one in the sense of consistency (same order of precision $\OO(\mu^2)$), where we remove the time dependency of the operator $\mathfrak{T}$. The structure of the time-independent model leads to a small improvement in terms of computational time due to its simple one-dimensional structure. In fact, this strategy was originally initiated for numerical simulations of the fully nonlinear and weakly dispersive GN models in the two-dimensional case~\cite{LannesMarche14}, in order to reduce significantly the computational time. \\
\\
We organize the paper as follows. 
In Section~\ref{FEsec}, we introduce the non-dimensonalized {\em full Euler system}. 
Section~\ref{GNCHsec} is devoted to recall the GN model in the Camassa-Holm regime, where an equivalent time-independent reformulation is given with an improved frequency dispersion based on the choice of a parameter $\alpha$. The stability issue of this new reformulation is discussed in the one and two layers cases. 
Section~\ref{NMSec} is then devoted to  the presentation of the numerical scheme. Firstly the hyperbolic/dispersive splitting is introduced. 
Then we describe the finite-volume spatial discretization of the hyperbolic part. The first-order finite volume scheme is given, an extension to the second-order accuracy is considered following the ``MUSCL" approach and finally high-order accuracy around discontinuities is achieved thanks to the implementation of a fifth-order WENO reconstruction. The finite-difference discretization of the dispersive part is then detailed and the boundary conditions are briefly described. Finally, we present in Section~\ref{NVsec} several numerical validations of our model. The one layer case is numerically validated: we compare the accuracy of the first, second and fifth orders of accuracy by studying the propagation of a solitary wave. To evaluate the influence of dispersive nonlinear terms, we consider the head-on collision of counter-propagating waves and the breaking of a Gaussian hump. Dealing with discontinuities is numerically validated by studying the dam-break problem in the one layer case. We highlight then the importance of the choice of the parameter $\alpha$ in improving the frequency dispersion of the model and compare our results with numerical experiments. Finally, we validate the ability of our numerical scheme to deal with discontinuities when considering a dam-break problem in the two layers case.  

\section{Full Euler system}\label{FEsec}
In this section, we briefly recall the derivation of the {\em full Euler system} governing the evolution equations of the two-layers flow and refer to 
 \cite{Anh09,BonaLannesSaut08,Duchene13,DucheneIsrawiTalhouk14} for more details. The two-layers flow considered are assumed to be incompressible, homogeneous, immiscible perfect fluids of different densities under the sole influence of gravity. 
\begin{figure}[H]
\centering 
\includegraphics[scale=1]{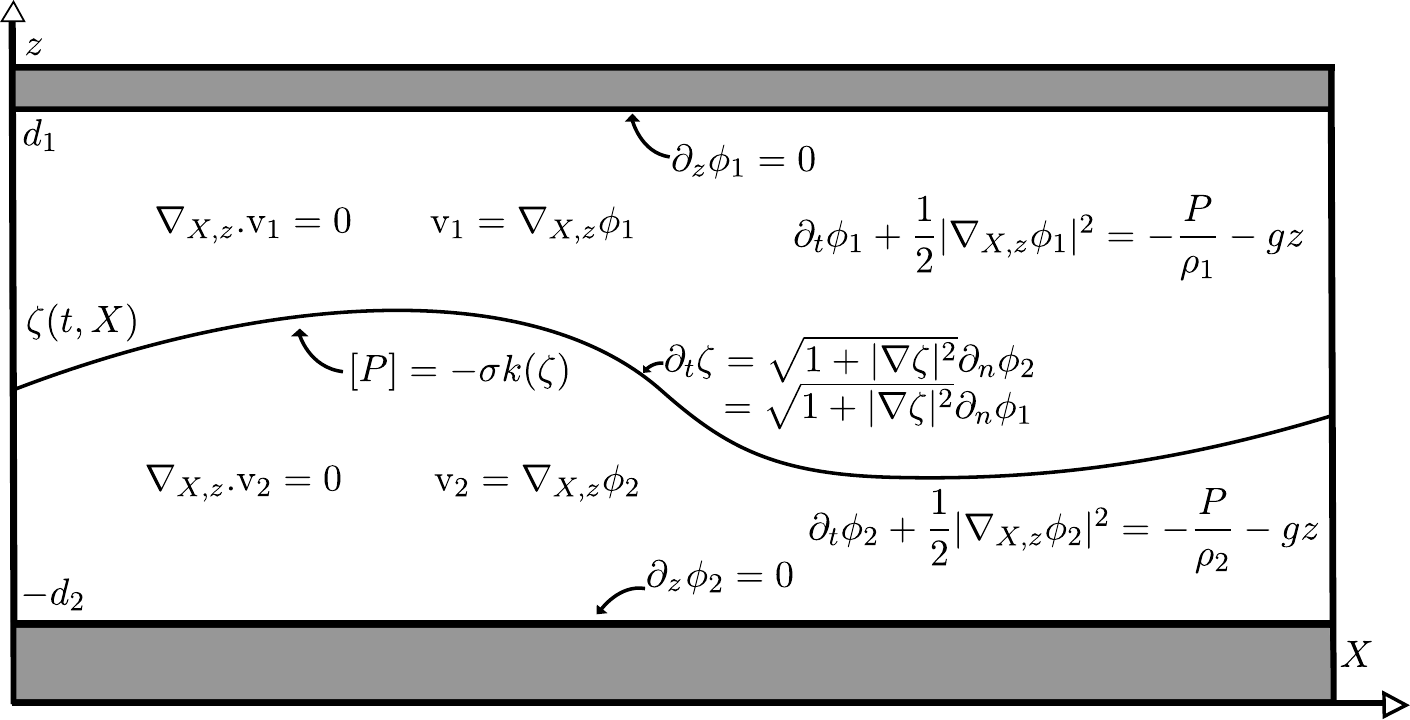}
\caption{Domain of study and governing equations.}
\label{domain}
\end{figure}
The study is restricted to the one-dimensional horizontal variable.
The deformation of the interface between the two layers is represented by the graph of function a $\zeta(t,x)$ while the bottom and the top surfaces are assumed to be rigid and flat. The domains of the upper and lower fluid at time $t$ (denoted, respectively, $\Omega_1^t$ and $\Omega_2^t$), are given by
\begin{eqnarray*}
 \Omega_1^t \ &=& \ \{\ (x,z)\in\RR\times\RR, \quad \zeta(t,x)\ \leq\ z\ \leq \ d_1\ \}, \\
 \Omega_2^t \ &=& \ \{\ (x,z)\in\RR\times\RR, \quad -d_2 \  \leq\ z\ \leq \ \zeta(t,x)\ \}.
\end{eqnarray*}
In what follows, we assume that the domains remain strictly connected, that is there exists $h_0>0$ such that $d_1-\zeta(t,x)\geq h_0>0$ and $d_2+\zeta(t,x)\geq h_0>0$.
 The density and velocity fields of the upper and lower layers are denoted by $(\rho_i,{\mathbf v}_i)$, $(i=1,2)$ respectively.
We assume the fluids to be incompressible, homogeneous and irrotational so that the velocity fields are divergence free and expressed as gradients of a potential denoted $\phi_i$. The fluids being ideal, that is with no viscosity, we may assume that each fluid follows the Euler equations. Assuming that the surface, the bottom or the interface are impenetrable one deduces the kinematic boundary conditions. Finally, the set of equations is completed by the continuity of the stress tensor at the interface. 

Altogether, the governing equations are given by the following system:
\everymath{\displaystyle}
  \begin{equation}  \label{eqn:EulerComplet}
\left\{\begin{array}{ll}
         \partial_x^2 \phi_i \ + \ \partial_z^2 \phi_i \ = \ 0 & \mbox{ in }\Omega^t_i, \ i=1,2,\\
         \partial_t \phi_i+\frac{1}{2} |\nabla_{x,z} \phi_i|^2=-\frac{P_{i}}{\rho_i}-gz & \mbox{ in }\Omega^t_i, \ i=1,2, \\
         \partial_{z}\phi_1 \ = \ 0  & \mbox{ on } \Gamma_{\rm t}\equiv\{(x,z),z=d_1\}, \\
         \partial_t \zeta  \ = \ \sqrt{1+|\partial_x\zeta|^2}\partial_{n}\phi_1 \ = \ \sqrt{1+|\partial_x\zeta|^2}\partial_{n}\phi_2  & \mbox{ on } \Gamma \equiv\{(x,z),z=\zeta(t,x)\},\\
         \partial_{z}\phi_2  \ = \ 0 &  \mbox{ on } \Gamma_{\rm b}\equiv\{(x,z),z=-d_2\}, \\
        \lim\limits_{\varepsilon\to 0} \Big(  P(t,x,\zeta(t,x)+\varepsilon)  
-  P(t,x,\zeta(t,x)-\varepsilon)  \Big)   = -\sigma k(\zeta) & \mbox{ on } \Gamma,
         \end{array}
\right.
\end{equation}
where $n$ denotes the unit upward normal vector at the interface.\\
The function $ k(\zeta)=-\partial_x \Big(\frac1{\sqrt{1+|\partial_x\zeta|^2}}\partial_x\zeta\Big)$ denotes the mean curvature of the interface and $\sigma$ the surface (or interfacial) tension coefficient.
\bigskip

To go further in the study of the flow of the two layers in order to build a numerical scheme for the asymptotic dynamics, we write the system in dimensionless form. To this end, we introduce dimensionless parameters and variables that reduce the setting to the physical regime under consideration.
Firstly, let $a$ be the maximum amplitude of the deformation of the interface and $\lambda$ the wavelength of the interface. Then the typical velocity of small propagating internal waves (or wave celerity) is given by
\[c_0 \ = \ \sqrt{g\frac{(\rho_2-\rho_1) d_1 d_2}{\rho_2 d_1+\rho_1 d_2}}.\]
Consequently, we introduce the dimensionless variables:
\[
 \t z \ \equiv\  \dfrac{z}{d_1}, \quad\quad  \t x\ \equiv \ \dfrac{x}{\lambda}, \quad\quad  \t t\ \equiv\ \dfrac{c_0}{\lambda}t,
\]
the dimensionless unknowns:
\[
 \t{\zeta}(\t t,\t x)\ \equiv\ \dfrac{\zeta(t,x)}{a}, \quad\quad \t{\phi_i}(\t t,\t x,\t z)\ \equiv\ \dfrac{d_1}{a\lambda c_0}\phi_i(t,x,z) \quad (i=1,2),
\]
and finally the dimensionless parameters:
\[
 \gamma\ =\ \dfrac{\rho_1}{\rho_2}, \quad \epsilon\ \equiv\ \dfrac{a}{d_1},\quad   \mu\ \equiv\ \dfrac{d_1^2}{\lambda^2},  \quad \delta\ \equiv \ \dfrac{d_1}{d_2},  \quad {\rm bo}\ =\ \frac{ g(\rho_2-\rho_1) d_1^2}{\sigma}.
\]
In what follows we assume that the depth ratio $\delta$ does not approach zero or infinity which means that the two layers of fluid have comparable depth. Therefore, the choice of the reference vertical length is harmless so we decided to choose $d_1$. We would like to mention that $\Bo=\mu\bo$ where $\Bo$ is the classical Bond number.
The significations of the different dimensionless parameters are the following :
\begin{itemize}
    \item $\gamma$ represents the ratio between the densities of the two layers,
    \item $\epsilon$ represents the amplitude of the interface between the two layers,
    \item $\mu$ represents the nonlinear effect of the shallowness approximation,
    \item $\delta$ represents the ratio between the depth of the two layers,
    \item $\bo$ represents the capillary effect of the interface.
\end{itemize}
Let us now remark that the system can be rewritten as two evolution equations coupling Zakharov's canonical variables \cite{Zakharov68,CraigSulem93}, $(\zeta,\psi)$ representing respectively the deformation of the interface and the trace of the dimensionless upper potential at the interface defined by $\psi \ \equiv \ \phi_1(t,x,\zeta(t,x)).$ \\
\\
In order to do so, we define the so-called Dirichlet-Neumann operators. The tildes are removed for the sake of readability.
\everymath{\displaystyle}
\begin{definition}[Dirichlet-Neumann operators]
Let $\zeta\in H^{t_0+1}(\RR)$, $t_0>1/2$, such that there exists $h>0$ with
$h_1 \ \equiv\  1-\epsilon\zeta \geq h>0$ and $h_2 \ \equiv \ \frac1\delta +\epsilon \zeta\geq h>0$, and let $\psi\in L^2_{\rm loc}(\RR),\partial_x \psi\in H^{1/2}(\RR)$.
 Then we define:
  \begin{eqnarray*}
 G^{\mu}\psi  &\equiv & G^{\mu}[\epsilon\zeta]\psi  \equiv  \sqrt{1+\mu|\epsilon\partial_x\zeta|^2}\big(\partial_n \phi_1 \big)\id{z=\epsilon\zeta}  =  -\mu\epsilon(\partial_x\zeta) (\partial_x\phi_1)\id{z=\epsilon\zeta}+(\partial_z\phi_1)\id{z=\epsilon\zeta},\\
H^{\mu,\delta}\psi  &\equiv  & H^{\mu,\delta}[\epsilon\zeta]\psi \equiv  \partial_x \big(\phi_2\id{z=\epsilon\zeta}\big)  =  (\partial_x\phi_2)\id{z=\epsilon\zeta}+\epsilon(\partial_x \zeta)(\partial_z\phi_2)\id{z=\epsilon\zeta},
\end{eqnarray*}
where, $\phi_1$ and $\phi_2$ are uniquely deduced from $(\zeta,\psi)$ as solutions of the following Laplace's problems:
\begin{eqnarray*}
&&\left\{
\begin{array}{ll}
 \left(\ \mu\partial_x^2 \ +\  \partial_z^2\ \right)\phi_1=0 & \mbox{ in } \Omega_1\equiv \{(x,z)\in \RR^{2},\ \epsilon{\zeta}(x)<z<1\}, \\
\partial_z \phi_1 =0  & \mbox{ on } \Gamma_{\rm t}\equiv \{(x,z)\in \RR^{2},\ z=1\}, \\
 \phi_1 =\psi & \mbox{ on } \Gamma\equiv \{(x,z)\in \RR^{2},\ z=\epsilon \zeta\},
\end{array}
\right.\\
&&\left\{
\begin{array}{ll}
 \left(\ \mu\partial_x^2\ + \ \partial_z^2\ \right)\phi_2=0 & \mbox{ in } \Omega_2\equiv\{(x,z)\in \RR^{2},\ -\frac{1}{\delta}<z<\epsilon\zeta\}, \\
\partial_{n}\phi_2 = \partial_{n}\phi_1 & \mbox{ on } \Gamma, \\
 \partial_{z}\phi_2 =0 & \mbox{ on } \Gamma_{\rm b}\equiv \{(x,z)\in \RR^{2},\ z=-\frac{1}{\delta}\}.
\end{array}
\right.
\end{eqnarray*}
\end{definition}
At this stage of the model, using the above definition and withouth making any assumption on  the different dimensionless parameters,
one can rewrite the nondimensionalized version of~\eqref{eqn:EulerComplet} as a system of two coupled evolution equations, which writes:
\begin{equation}\label{eqn:EulerCompletAdim}
\left\{ \begin{array}{l}
\displaystyle\partial_{ t}{\zeta} \ -\ \frac{1}{\mu}G^{\mu}\psi\ =\ 0,  \\ \\
\displaystyle\partial_{ t}\Big(H^{\mu,\delta}\psi-\gamma \partial_x{\psi} \Big)\ + \ (\gamma+\delta)\partial_x{\zeta} \ + \ \frac{\epsilon}{2} \partial_x\Big(|H^{\mu,\delta}\psi|^2 -\gamma |\partial_x {\psi}|^2 \Big) \\
\displaystyle\hspace{5cm} = \ \mu\epsilon\partial_x\N^{\mu,\delta}-\mu\frac{\gamma+\delta}{\bo}\frac{\partial_x \big(k(\epsilon\sqrt\mu\zeta)\big)}{{\epsilon\sqrt\mu}} \ ,
\end{array}
\right.
\end{equation}
where we denote:
\[  \N^{\mu,\delta} \ \equiv \ \frac{\big(\frac{1}{\mu}G^{\mu}\psi+\epsilon(\partial_x{\zeta})H^{\mu,\delta}\psi \big)^2\ -\ \gamma\big(\frac{1}{\mu}G^{\mu}\psi+\epsilon(\partial_x{\zeta})(\partial_x{\psi}) \big)^2}{2(1+\mu|\epsilon\partial_x{\zeta}|^2)}.
      \]
We will refer to \eqref{eqn:EulerCompletAdim} as the {\em full Euler system}.

\section{Green-Naghdi model in the Camassa-Holm regime}\label{GNCHsec}
We now recall the new Green-Naghdi model in the Camassa-Holm regime recently derived by Duchêne, Israwi and Talhouk in \cite{DucheneIsrawiTalhouk15}.
This new model is derived after expanding the different operators of the original Green-Naghdi model with respect to $\epsilon \mbox{ and } \mu$. Then several additional transformations are made using the smallness assumption $\epsilon=\OO(\sqrt{\mu})$ in order to produce an equivalent precise system of partial differential equations whose unknowns are
$\zeta \mbox{  and }v$ :

\everymath{\displaystyle}
\begin{equation}\label{eq:Serre2mr}\left\{ \begin{array}{l}
\partial_{ t}\zeta +\partial_x\left(\dfrac{h_1 h_2}{h_1+\gamma h_2}v\right)\ =\  0,\\ \\
\mathfrak T[\epsilon\zeta] \left( \partial_{ t} v + \epsilon \varsigma {v } \partial_x {v} \right) + (\gamma+\delta)q_1(\epsilon\zeta)\partial_x \zeta   \\
\qquad \qquad+\frac\epsilon2 q_1(\epsilon\zeta) \partial_x \left(\frac{h_1^2  -\gamma h_2^2 }{(h_1+\gamma h_2)^2}| v|^2-\varsigma |v|^2\right)=  -  \mu \epsilon\frac23\frac{1-\gamma}{(\gamma+\delta)^2} \partial_x\big((\partial_x v)^2\big) ,
\end{array} \right. \end{equation}
where $h_1\equiv1-\epsilon\zeta$ (resp. $h_2\equiv\frac1\delta+\epsilon\zeta$) denotes the depth of the upper (resp. lower) fluid, and $v$ is the {\em shear mean velocity} defined by:
 \begin{equation*}v\equiv \frac{1}{h_2}\int_{-\frac1\delta}^{\epsilon\zeta(t,x)} \partial_x \phi_2(t,x,z) \ dz - \frac{\gamma}{h_1}\int_{\epsilon\zeta(t,x)}^{1} \partial_x \phi_1(t,x,z) \ dz.\end{equation*}
The operator ${\mathfrak T}$ is defined as:
\begin{equation*}
    {\mathfrak T}[\epsilon\zeta]V \ = \ q_1(\epsilon\zeta)V \ - \ \mu\nu \partial_x \Big(q_2(\epsilon\zeta)\partial_xV \Big),
\end{equation*}
with $q_i(X)\equiv 1+\kappa_i X $ ($i=1,2$) and $\nu$, $\kappa_1$, $\kappa_2$, $\varsigma$ are defined  below.
Let us first introduce the following constants in order to ease the reading:
\begin{equation}\label{eqn:deflambdaalphabeta}
    \lambda=\frac{1+\gamma\delta}{3\delta(\gamma+\delta)}\ , \quad 
    \alpha=\dfrac{1-\gamma}{(\gamma+\delta)^2} \quad 
    \mbox{ and } \quad \beta=\dfrac{(1+\gamma\delta)(\delta^2-\gamma)}{\delta(\gamma+\delta)^3}\ .
\end{equation}
Thus
\begin{equation*}
    \nu \ = \ \lambda-\frac1{\bo} \ = \ \frac{1+\gamma\delta}{3\delta(\gamma+\delta)}-\frac1{\bo},
\end{equation*}
\begin{equation*}
    (\lambda-\frac1{\bo})\kappa_1 \ = \ \frac{\gamma+\delta}{3}(2\beta-\alpha) , \quad   (\lambda-\frac1{\bo})\kappa_2 \ = \ (\gamma+\delta)\beta ,    \end{equation*}
\begin{equation*}
    (\lambda-\frac1{\bo})\varsigma \ = \ \frac{2\alpha-\beta}{3} \ - \ \frac{1}{\bo} \dfrac{\delta^2 -\gamma }{(\delta+\gamma )^2}.
\end{equation*}
The function $f$ is defined as follows:
\[ f:X\to\frac{(1-X)(\delta^{-1}+X)}{1-X+\gamma(\delta^{-1}+X)},\]
so one has:
\[
f(\epsilon\zeta) \ =\ \dsp\frac{h_1h_2}{h_1+\gamma h_2} \quad \mbox{ and } \quad f'(\epsilon\zeta) \ =\ \dsp\frac{h_1^2-\gamma h_2^2}{(h_1+\gamma h_2)^2} \ .\]
Additionally, let us denote:
\begin{equation}\label{defkappaq3}
\kappa=\frac23\frac{1-\gamma}{(\delta+\gamma)^2}
\quad \mbox{ and } \quad q_3(\epsilon\zeta)=\frac12\big(f'(\epsilon\zeta)-\varsigma\big) = \frac12\Big(\frac{h_1^2-\gamma h_2^2}{(h_1+\gamma h_2)^2}-\varsigma\Big) \ .
\end{equation}
Problem \eqref{eq:Serre2mr} writes now in compact form:
\begin{equation}\label{GNCH2}\left\{ \begin{array}{l}
\dsp \partial_{ t}\zeta +\partial_x\big(f(\epsilon\zeta)  v\big)\ =\  0,\\ \\
\dsp {\mathfrak T} \left( \partial_{ t}  v + \epsilon\varsigma v \partial_x v  \right) + (\gamma+\delta)q_1(\epsilon\zeta)\partial_x \zeta + \epsilon q_1(\epsilon\zeta)\partial_x(q_3(\epsilon\zeta)  {v}^2)  + \mu\epsilon\kappa\partial_x\big((\partial_x v)^2\big)=0.
\end{array} \right. \end{equation}
Let us now recall the regime of validity of the system~\eqref{GNCH2} as exactly given in~\cite{DucheneIsrawiTalhouk15}. In the first place, we consider the so-called {\em shallow water regime} for two layers of comparable depths:
\begin{multline} \label{eqn:defRegimeSWmr}
\P_{\rm SW} \ \equiv \ \Big\{ (\mu,\epsilon,\delta,\gamma,\bo):\ 0\ < \ \mu \ \leq \ \mu_{\max}, \ 0 \ \leq \ \epsilon \ \leq \ 1, \ \delta \in (\delta_{\min},\delta_{\max}),  \Big.\\
 \Big. \ 0\ \leq \ \gamma\ <\ 1,\  \bo_{\min}\leq \bo\leq \infty \ \Big\},
\end{multline}
with given $0\leq \mu_{\max},\delta_{\min}^{-1},\delta_{\max},\bo_{\min}^{-1}<\infty$. 

Two additional key restrictions are necessary for the validity of the model:
\begin{equation} \label{eqn:defRegimeCHmr}
\P_{\rm CH}  \equiv  \left\{ (\mu,\epsilon,\delta,\gamma,\bo) \in \P_{\rm SW}:\ \epsilon  \leq  M \sqrt{\mu}  \quad \mbox{ and } \quad \nu \equiv  \frac{1+\gamma\delta}{3\delta(\gamma+\delta)}-\frac1{\bo} \ge  \nu_{0}  \ \right\},
\end{equation}
with given $0\leq M,\nu_0^{-1}<\infty$.\\
\\
Duch\^ene, Israwi and Talhouk proved in~\cite{DucheneIsrawiTalhouk15}, that the system \eqref{GNCH2} is well-posed (in the sense of Hadamard) in the energy space $X^s\ = \ H^s(\RR)\times H^{s+1}(\RR)$, endowed with the norm
\[
\forall\; U=(\zeta,v)^\top \in X^s, \quad \vert U\vert^2_{X^s}\equiv \vert \zeta\vert^2 _{H^s}+\vert v\vert^2 _{H^s}+ \mu\vert \partial_xv\vert^2 _{H^s}.
\]
This result is obtained for parameters in the {\em Camassa-Holm regime}~\eqref{eqn:defRegimeCHmr} and under the following non-zero depth and ellipticity (for the operator $\mathfrak{T}$) conditions:
\begin{equation}\tag{H1}
\exists h_{01}>0 \mbox{ such that, } \inf_{x\in \RR} h_1 \geq\ h_{01}\ >\ 0, \quad \inf_{x\in \RR} h_2\geq\ h_{01}\ >\ 0.
\end{equation}
\begin{equation}\tag{H2}
\exists h_{02}>0 \mbox{ such that, } \inf_{x\in \RR} \left(1+\epsilon\kappa_1\zeta\right) \ge \ h_{02} \ > \ 0, \  \quad    \inf_{x\in \RR}  \left( 1+\epsilon \kappa_2\zeta  \right)\ge \ h_{02} \ > \ 0.
\end{equation}
They also prove that this new asymptotic model model is fully justified by a convergence result in the Camassa-Holm regime.

\subsection{Reformulation of the model}\label{REFsec}
One can easily check that the operator $\mathfrak{T}$ can be written under the form:
\begin{equation*}
\mfT[\epsilon\zeta]V=(q_1(\epsilon\zeta)I+\mu\nu T[\epsilon\zeta])V \quad \mbox{with} \quad T[\epsilon\zeta]V=-\partial_x(q_2(\epsilon\zeta)\partial_x V).
\end{equation*}
Therefore, this operator is time dependant through the presence of the term $\zeta$ and at each time, this operator has to be inverted in order to 
solve equation \eqref{GNCH2}
Thus we want to derive a new model equivalent to \eqref{eq:Serre2mr} up to a lower order in $\epsilon \mbox{ or } \mu$ but with a structure 
adapted to construct a numerical scheme for its resolution.

Let us first denote : 
\begin{equation}
\label{Q1}Q_1(v)=\kappa\partial_x((\partial_xv)^2),
\end{equation}
such that one can rewrite problem  \eqref{GNCH2}  as : 
\begin{equation}\label{GNCH3}
\left\{ \begin{array}{l}
\hspace*{-0.25cm}\dsp \partial_{ t}\zeta +\partial_x\big(f(\epsilon\zeta)  v\big)\ =\  0,\\ \\
\hspace*{-0.25cm}\dsp \Big(q_1(\epsilon\zeta)I+\mu\nu T[\epsilon\zeta]\Big) \Big( \partial_{ t}  v + \epsilon\varsigma v \partial_x v  \Big) + (\gamma+\delta)q_1(\epsilon\zeta)\partial_x \zeta + \epsilon q_1(\epsilon\zeta)\partial_x(q_3(\epsilon\zeta)  {v}^2)  + \mu\epsilon Q_1(v)=0 .
\end{array} \right. \end{equation}
As we said before since the operator $(q_1(\epsilon\zeta)I+\mu\nu T[\epsilon\zeta])$ depends on $\zeta(t,x)$, 
one has to remove this time dependency to avoid its inversion at each time.

Firstly, let us write the operator $T[\epsilon\zeta]$ under the form:
\[T[\epsilon\zeta]V=-\partial_x (q_2(\epsilon\zeta)\partial_x V)=T[0]V+ \epsilon S[\zeta]V,\]with
\[
T[0]V=-\partial_x^2V
\quad \mbox{and} \quad S[\zeta]V=-\kappa_2\partial_x(\zeta\partial_xV).\]
Expanding the term $\Big(q_1(\epsilon\zeta)I+\mu\nu T[\epsilon\zeta]\Big) \Big( \partial_{ t}  v + \epsilon\varsigma v \partial_x v  \Big)$ in problem
\eqref{GNCH3} leads to :
\begin{equation}\label{GNCH4}\left\{ \begin{array}{l}
\dsp \partial_{ t}\zeta +\partial_x\big(f(\epsilon\zeta)  v\big)\ =\  0,\\ \\
\dsp \Big(I+\mu\nu T[0]\Big) \Big( \partial_{ t}  v + \epsilon\varsigma v \partial_x v  \Big) + (\gamma+\delta)q_1(\epsilon\zeta)\partial_x \zeta + \epsilon q_1(\epsilon\zeta)\partial_x(q_3(\epsilon\zeta)  {v}^2)  + \mu\epsilon Q_1(v)\\+\mu\epsilon \nu  S[\zeta] \left( \partial_{ t}  v + \epsilon\varsigma v \partial_x v  \right)+\epsilon\kappa_1\zeta \left( \partial_{ t}  v + \epsilon\varsigma v \partial_x v  \right)=0 .
\end{array} \right. \end{equation}
But we have: 
$\mu\epsilon \nu  S[\zeta] \left( \partial_{ t}  v + \epsilon\varsigma v \partial_x v  \right) = \mu\epsilon \nu S[\zeta] \partial_{ t}  v + \OO(\mu\epsilon^2)$.
By assumption the term $\OO(\mu\epsilon^2)$ is of order $\OO(\mu^2)$; thus Problem \eqref{GNCH4} is equivalent to the following system:
\begin{equation}\label{GNCH5}\left\{ \begin{array}{l}
\dsp \partial_{ t}\zeta +\partial_x\big(f(\epsilon\zeta)  v\big)\ =\  0,\\ \\
\dsp \Big(I+\mu\nu T[0]\Big) \Big( \partial_{ t}  v + \epsilon\varsigma v \partial_x v  \Big) + (\gamma+\delta)q_1(\epsilon\zeta)\partial_x \zeta + \epsilon q_1(\epsilon\zeta)\partial_x(q_3(\epsilon\zeta)  {v}^2)  + \mu\epsilon Q_1(v)\\+\mu\epsilon \nu  S[\zeta] \left( \partial_{ t}  v  \right)+\epsilon\kappa_1\zeta \left( \partial_{ t}  v + \epsilon\varsigma v \partial_x v  \right)=0 .
\end{array} \right. \end{equation}

\subsection{Improved Green-Naghdi equations}\label{GNID}
As done by many authors, see \cite{Witting84,MMS91,CBB07}, the frequency dispersion of problem \eqref{GNCH5} can be improved by adding some terms of order $\OO(\mu^2)$ to the momentum equation. The accuracy of the model is not affected by this manipulation since this equation is precise up to terms of the same order as the added ones, namely $\OO(\mu^2)$.

Going back to problem \eqref{GNCH3},  one notices that:
\begin{equation}\label{eq1}q_1(\epsilon\zeta)\left( \partial_{ t}  v + \epsilon\varsigma v \partial_x v  \right)+(\gamma+\delta)q_1(\epsilon\zeta)\partial_x \zeta + \epsilon q_1(\epsilon\zeta)\partial_x(q_3(\epsilon\zeta)  {v}^2)+\OO(\mu)=0.\end{equation}
We can now divide \eqref{eq1} by $q_1(\epsilon\zeta)$ to obtain:
\begin{equation}\label{eq2}\partial_{ t}  v=-\epsilon\varsigma v \partial_x v -(\gamma+\delta)\partial_x \zeta -\epsilon\partial_x(q_3(\epsilon\zeta)  {v}^2)-\dfrac{\OO(\mu)}{q_1(\epsilon\zeta)}.\end{equation}
Let us fix $\alpha\in\RR$. Multiplying \eqref{eq2} by $(1-\alpha)$ leads to:
\begin{equation}\label{eq3}\partial_{ t}  v=\alpha\partial_{ t}  v-(1-\alpha)[\epsilon\varsigma v \partial_x v -(\gamma+\delta)\partial_x \zeta -\epsilon\partial_x(q_3(\epsilon\zeta)  {v}^2)]-(1-\alpha)\dfrac{\OO(\mu)}{q_1(\epsilon\zeta)}.\end{equation}
Replacing $\partial_t v$ by its expression given in \eqref{eq3} in the second equation of \eqref{GNCH5} and regrouping the $\OO(\mu^2)$ 
terms yields to the following equation:
\begin{multline*}
(I+\mu\nu T[0]) \left( \alpha\partial_{ t}  v-(1-\alpha)[\epsilon\varsigma v \partial_x v +(\gamma+\delta)\partial_x \zeta +\epsilon\partial_x(q_3(\epsilon\zeta)  {v}^2)]-(1-\alpha)\dfrac{\OO(\mu)}{q_1(\epsilon\zeta)} + \epsilon\varsigma v \partial_x v  \right) \\ + (\gamma+\delta)q_1(\epsilon\zeta)\partial_x \zeta + \epsilon q_1(\epsilon\zeta)\partial_x(q_3(\epsilon\zeta)  {v}^2)  + \mu\epsilon Q_1(v)+\mu\epsilon\nu S[\zeta] \partial_t v+\epsilon\kappa_1\zeta \left( \partial_{ t}  v + \epsilon\varsigma v \partial_x v \right)+ \OO(\mu^2)=0
\end{multline*}
After straightforward computations,
\begin{multline}\label{eq4}
(I+\mu\nu\alpha T[0]) \left( \partial_{ t}  v+\epsilon\varsigma v\partial_x v\right) +(\alpha-1)[\partial_t v+\epsilon\varsigma v \partial_x v +(\gamma+\delta)\partial_x \zeta +\epsilon\partial_x(q_3(\epsilon\zeta)  {v}^2)]\\+(\alpha-1)\mu\nu T[0]\big((\gamma+\delta)\partial_x \zeta +\epsilon\partial_x(q_3(\epsilon\zeta)  {v}^2)\big)+(\gamma+\delta)q_1(\epsilon\zeta)\partial_x \zeta+\epsilon q_1(\epsilon\zeta)\partial_x(q_3(\epsilon\zeta)v^2)\\+\mu\epsilon Q_1(v)+ \mu \epsilon \nu S[\zeta] \partial_t v+\epsilon\kappa_1\zeta \left( \partial_{ t}  v + \epsilon\varsigma v \partial_x v  \right)-(1-\alpha)\dfrac{\OO(\mu)}{q_1(\epsilon\zeta)}+\OO(\mu^2) =0
\end{multline}
Replacing in \eqref{eq4},
$\epsilon\kappa_1\zeta \left( \partial_{ t}  v + \epsilon\varsigma v \partial_x v  \right)$ by $\epsilon\kappa_1\zeta \Big(-(\gamma+\delta)\partial_x \zeta -\epsilon\partial_x(q_3(\epsilon\zeta)  {v}^2)-\dfrac{\OO(\mu)}{q_1(\epsilon\zeta)}\Big)$,
we have:
\begin{multline}\label{eq5}
(I+\mu\nu\alpha T[0]) \left( \partial_{ t}  v+\epsilon\varsigma v\partial_x v\right) +(I-\mu\nu(1-\alpha) T[0])(\gamma+\delta)\partial_x \zeta +\epsilon\partial_x(q_3(\epsilon\zeta)  {v}^2))\\+\mu\epsilon Q_1(v)+
\mu\epsilon \nu S[\zeta] \partial_t v-\epsilon\kappa_1\zeta\dfrac{\OO(\mu)}{q_1(\epsilon\zeta)}+\OO(\mu^2) =0.
\end{multline}
One deduces easily that,
\begin{multline}\label{eq6}(I-\mu\nu(1-\alpha) T[0])(\gamma+\delta)\partial_x \zeta +\epsilon\partial_x(q_3(\epsilon\zeta)  {v}^2))=\dfrac{1}{\alpha}\big((\gamma+\delta)\partial_x \zeta +\epsilon\partial_x(q_3(\epsilon\zeta)  {v}^2)\big) \\ +\dfrac{\alpha-1}{\alpha}(I+\mu\nu\alpha T[0])((\gamma+\delta)\partial_x \zeta +\epsilon\partial_x(q_3(\epsilon\zeta)  {v}^2)).
\end{multline}
Plugging \eqref{eq6}  in \eqref{eq5} one has:
\begin{multline}\label{eq7}
(I+\mu\nu\alpha T[0]) \big[ \partial_{ t}  v+\epsilon\varsigma v\partial_x v+\dfrac{\alpha-1}{\alpha}\big((\gamma+\delta)\partial_x \zeta +\epsilon\partial_x(q_3(\epsilon\zeta)  {v}^2)\big)\big]+\dfrac{1}{\alpha}\big((\gamma+\delta)\partial_x \zeta +\epsilon\partial_x(q_3(\epsilon\zeta)  {v}^2)\big)\\+\mu\epsilon Q_1(v)+\mu \epsilon \nu S[\zeta] \partial_t v-\epsilon\kappa_1\zeta\dfrac{\OO(\mu)}{q_1(\epsilon\zeta)}+\OO(\mu^2) =0
\end{multline}
Remembering  that $\mu\nu T[0] ( \partial_{ t}  v+\epsilon\varsigma v\partial_x v) +\mu\epsilon Q_1(v)+\mu\epsilon \nu S[\zeta] \partial_t v= \OO(\mu)$, 
and regrouping the $\OO(\mu\epsilon^2)$ terms, \eqref{eq7} becomes:
\begin{multline}\label{eq8}
(I+\mu\nu\alpha T[0]) \big[ \partial_{ t}  v+\epsilon\varsigma v\partial_x v+\dfrac{\alpha-1}{\alpha}\big((\gamma+\delta)\partial_x \zeta +\epsilon\partial_x(q_3(\epsilon\zeta)  {v}^2)\big)\big]+\dfrac{1}{\alpha}\big((\gamma+\delta)\partial_x \zeta +\epsilon\partial_x(q_3(\epsilon\zeta)  {v}^2)\big)\\+\mu\epsilon Q_1(v)+\mu \epsilon \nu  S[\zeta] (\partial_t v)-\epsilon\kappa_1\zeta\mu\nu T[0]
(\partial_t v)+\OO(\mu^2,\mu\epsilon^2) =0
\end{multline}
Expanding in terms of $\epsilon$ and $\mu$,  and looking for the term $\partial_t v$, we get:
 \begin{equation}\label{dtv1}
 \partial_{ t}  v = -(I+\mu\nu \alpha T[0])^{-1}[(\gamma+\delta)\partial_x \zeta]+\OO(\epsilon,\mu)
 \end{equation}
Replacing $ \partial_{ t}  v $ by its expression obtained in~\eqref{dtv1} in the last two terms of the equation \eqref{eq8} the Green–Naghdi equations with improved dispersion can therefore be written as (up to the order $\OO(\mu^2)$) :
\begin{equation}\label{GNCH6}\left\{ \begin{array}{l}
\dsp \partial_{ t}\zeta +\partial_x\big(f(\epsilon\zeta)  v\big)\ =\  0,\\ \\
\dsp (I+\mu\nu\alpha T[0]) \big[ \partial_{ t}  v+\epsilon\varsigma v\partial_x v+\dfrac{\alpha-1}{\alpha}\big((\gamma+\delta)\partial_x \zeta +\epsilon\partial_x(q_3(\epsilon\zeta)  {v}^2)\big)\big]\\+\dfrac{1}{\alpha}\big((\gamma+\delta)\partial_x \zeta +\epsilon\partial_x(q_3(\epsilon\zeta)  {v}^2)\big)+\mu\epsilon Q_1(v)+\mu\epsilon\nu Q_2(\zeta)+\mu\epsilon\nu Q_3(\zeta) =0.
\end{array} \right. \end{equation}
with 
\begin{equation}\label{Q2}
Q_2(\zeta)=- S[\zeta] \Big(I+\mu\nu \alpha T[0]\Big)^{-1}[(\gamma+\delta)\partial_x \zeta ] ,
\end{equation}
and
\begin{equation}\label{Q3}
Q_3(\zeta)=\kappa_1\zeta T[0]\Big(I+\mu\nu \alpha T[0]\Big)^{-1}[(\gamma+\delta)\partial_x \zeta ].
\end{equation}
In \cite{CLM}, a three-parameter family of Green–Naghdi equations in the one layer case is derived yielding additional improvements of the dispersive properties. In this paper we stick to the one-parameter family~\eqref{GNCH6}. This new formulation does not contain any third-order derivative, thus one can expect more stable and robust numerical computations.

\subsection{Choice of the parameter $\alpha$}\label{Secalphachoice}
The choice of the parameter $\alpha$ is motivated by the agreement of the dispersion properties of the \emph {full Euler system} and the improved Green-Naghdi system
\eqref{GNCH6} in term of the dispersion relation. Therefore, we have first to find the dispersion relation for the improved Green-Naghdi system, then the dispersion relation
for the \emph {full Euler system} and finally to find an optimal parameter $\alpha_{opt}$ to ensure a good agreement between these two relations.
\subsubsection{The dispersion relation associated to the improved GN formulation}
The system \eqref{GNCH6} can be written under the following form:
\begin{equation}\label{GNCH7}\left\{ \begin{array}{l}
\dsp \partial_{ t}\zeta +\partial_x\big(f(\epsilon\zeta)  v\big)\ =\  0,\\ \\
\dsp (I+\mu\nu\alpha T[0]) \big[ \partial_{ t}  v+\epsilon\varsigma v\partial_x v+\dfrac{\alpha-1}{\alpha}\big((\gamma+\delta)\partial_x \zeta +\epsilon\partial_x(q_3(\epsilon\zeta)  {v}^2)\big)\big]\\+\dfrac{1}{\alpha}\big((\gamma+\delta)\partial_x \zeta +\epsilon\partial_x(q_3(\epsilon\zeta)  {v}^2)\big)+\mu\epsilon Q_1(v)+\mu\epsilon\nu Q_2(\zeta)+\mu\epsilon\nu Q_3(\zeta) =0,
\end{array} \right. \end{equation}
where the operators $Q_1$, $Q_2$ and $Q_3$ are explicitly given in~\eqref{Q1},~\eqref{Q2} and~\eqref{Q3}. Looking at the linearization of~\eqref{GNCH7} around the rest state $(\zeta,v)=(0,0)$, one derives the dispersion
relation associated to~\eqref{GNCH7}.\\
This relation is obtained by looking for plane wave solutions of the form 
$(\zeta,v) = (\zetabar,\vbar)e^{i(kx-wt)}$, with $k$ the spatial wave number and $\omega$ the time pulsation, of  the linearized equations.\\
\\
The first equation of \eqref{GNCH7} gives:
\[-iw\zetabar+\dfrac{1}{(\gamma+\delta)}ik\vbar=0 \ .\]
Thus we first obtain:
\begin{equation}\label{eq9}\vbar=\dfrac{w(\gamma+\delta)}{k}\zetabar.\end{equation}
From the second equation of \eqref{GNCH7} we have:
\[(1+\mu\nu\alpha k^2)\Big(-iw\vbar+\dfrac{\alpha-1}{\alpha}(\gamma+\delta)ik\zetabar\Big)+\dfrac{1}{\alpha}(\gamma+\delta)ik\zetabar=0 .\]
This equation may be written as:
\begin{equation}\label{eq10}
-iw\vbar+(\gamma+\delta)ik\zetabar-\mu\nu\alpha iwk^2\vbar+\mu\nu(\alpha-1)(\gamma+\delta)ik^3\zetabar=0.
\end{equation}
Replacing $\vbar$ in~\eqref{eq10} by its expression given in~\eqref{eq9} we obtain:
\begin{equation}
-w^2\dfrac{(\gamma+\delta)}{k}\zetabar+(\gamma+\delta)k\zetabar-\mu\nu\alpha wk^2\dfrac{w(\gamma+\delta)}{k}\zetabar+\mu\nu(\alpha-1)(\gamma+\delta)k^3\zetabar=0.
\end{equation}
After straightforward computations, we obtain the following dispersion relation of the new GN formulation for $\alpha\geq 1$:
\begin{equation}\label{rdgn}
w_{\alpha,GN}=\pm|k|\sqrt{\dfrac{1+\mu\nu(\alpha-1)k^2}{1+\mu\nu\alpha k^2}}
\end{equation}
Defining the linear phase velocity associated to~\eqref{rdgn} as:
\[C_{GN}^p(k)=\dfrac{w(k)}{|k|} ,\]
we choose $\alpha$ such that the phase velocity stays close to the reference phase velocity $C_S^p(k)$ coming from Stokes
linear theory. Thanks to this approach, the error on the phase velocity is minimized for any discrete value of $\mu|k|$ and the corresponding local optimal
value of $\alpha$, denoted by $\alpha_{opt}$ is computed.
In order to do so, we have firstly to obtain the dispersion relation of the original \emph {full Euler system}.
\subsubsection{The dispersion relation associated to the full-Euler system}
Let us recall the {\em full-Euler system} given in~\eqref{eqn:EulerCompletAdim}:
\begin{equation}\label{Eulercomplet}
\left\{ \begin{array}{l}
\displaystyle\partial_{ t}{\zeta} \ -\ \frac{1}{\mu}G^{\mu}\psi\ =\ 0,  \\ \\
\displaystyle\partial_{ t}\Big(H^{\mu,\delta}\psi-\gamma \partial_x{\psi} \Big)\ + \ (\gamma+\delta)\partial_x{\zeta} \ + \ \frac{\epsilon}{2} \partial_x\Big(|H^{\mu,\delta}\psi|^2 -\gamma |\partial_x {\psi}|^2 \Big) \\ \hspace{5cm} = \mu\epsilon\partial_x\N^{\mu,\delta}-\mu\frac{\gamma+\delta}{\bo}\frac{\partial_x \big(k(\epsilon\sqrt\mu\zeta)\big)}{{\epsilon\sqrt\mu}} \ ,
\end{array}
\right.
\end{equation}
where
\[  \N^{\mu,\delta} \ \equiv \ \frac{\big(\frac{1}{\mu}G^{\mu}\psi+\epsilon(\partial_x{\zeta})H^{\mu,\delta}\psi \big)^2\ -\ \gamma\big(\frac{1}{\mu}G^{\mu}\psi+\epsilon(\partial_x{\zeta})(\partial_x{\psi}) \big)^2}{2(1+\mu|\epsilon\partial_x{\zeta}|^2)},
      \]
and
\begin{eqnarray*}
G^{\mu}\psi  &\equiv & G^{\mu}[\epsilon\zeta]\psi  \equiv  \sqrt{1+\mu|\epsilon\partial_x\zeta|^2}\big(\partial_n \phi_1 \big)\id{z=\epsilon\zeta}  =  -\mu\epsilon(\partial_x\zeta) (\partial_x\phi_1)\id{z=\epsilon\zeta}+(\partial_z\phi_1)\id{z=\epsilon\zeta},\\
H^{\mu,\delta}\psi  &\equiv  & H^{\mu,\delta}[\epsilon\zeta,\beta b]\psi \equiv  \partial_x \big(\phi_2\id{z=\epsilon\zeta}\big)  =  (\partial_x\phi_2)\id{z=\epsilon\zeta}+\epsilon(\partial_x \zeta)(\partial_z\phi_2)\id{z=\epsilon\zeta},
\end{eqnarray*}
where $\phi_1$ and $\phi_2$ are uniquely defined (up to a constant for $\phi_2$) as the solutions in $H^2(\RR)$ of the  following Laplace's problems.

\begin{eqnarray}
\label{Laplace1} &&\left \{
\begin{array}{ll}\left(\ \mu\partial_x^2 \ +\  \partial_z^2\ \right)\ \phi_1=0 & \mbox{ in } \Omega_1\equiv \{(x,z)\in \RR^{2},\ \epsilon{\zeta}(x)<z<1\}, \\
\partial_z \phi_1 =0  & \mbox{ on } \Gamma_{\rm t}\equiv \{(x,z)\in \RR^{2},\ z=1\}, \\
\phi_1 =\psi & \mbox{ on } \Gamma\equiv \{(x,z)\in \RR^{2},\ z=\epsilon \zeta\},
\end{array}
\right.\\
\label{Laplace2}&&\left\{
\begin{array}{ll}
\left(\ \mu\partial_x^2\ + \ \partial_z^2\ \right)\ \phi_2=0 & \mbox{ in } \Omega_2\equiv\{(x,z)\in \RR^{2},\ -\frac{1}{\delta}<z<\epsilon\zeta\}, \\
\partial_{n}\phi_2 = \partial_{n}\phi_1 & \mbox{ on } \Gamma,\\
\partial_{z}\phi_2 =0 & \mbox{ on } \Gamma_{\rm b}\equiv \{(x,z)\in \RR^{2},\ z=-\frac{1}{\delta}\}.
\end{array}
\right.
\end{eqnarray}
To obtain the  dispersion relation associated to the {\em full-Euler system}, we first have to linearize it around the rest state that is
$\zeta=0 \, , \, v = 0$. We have thus to compute the two operators  $G^{\mu}$ and $H^{\mu,\delta}$ at the rest state. 
Therefore we have firstly to linearize the two previous equations \eqref{Laplace1}-\eqref{Laplace2}, 
then write the linear system in the wave number space by performing a Fourier transform. 
After performing some algebraic computations, we will obtain the dispersion relation. 
Let us remark that the two equations \eqref{Laplace1} and \eqref{Laplace2} are very similar except for the boundary conditions,
and the space domain $\Omega_1$ and $\Omega_2$. Therefore we have chosen to perform the complete analysis for the two equations.

$\bullet$ Linearizing~\eqref{Laplace1} around the rest state ($\zeta=0$) gives:
\begin{eqnarray}
\label{Laplace1lin} &\left\{
\begin{array}{ll}\left(\ \mu\partial_x^2 \ +\  \partial_z^2\ \right)\ \phi_1=0 & \mbox{ in } \Omega_1\equiv \{(x,z)\in \RR^{2},\ 0<z<1\}, \\
\partial_z \phi_1 =0  & \mbox{ on } \Gamma_{\rm t}\equiv \{(x,z)\in \RR^{2},\ z=1\}, \\
\phi_1 =\psi & \mbox{ on } \Gamma\equiv \{(x,z)\in \RR^{2},\ z=0\},
\end{array}
\right.
\end{eqnarray}
Applying the Fourier transform with respect to $x$, one has:
\begin{eqnarray}
\label{Laplace1four} &\left\{
\begin{array}{ll}\ -\mu|k|^2\widehat{\phi_1}+\partial_z^2\widehat{\phi_1}=0 & \mbox{ in } \Omega_1\equiv \{(x,z)\in \RR^{2},\ 0<z<1\}, \\
\partial_z \widehat{\phi_1} =0  & \mbox{ on } \Gamma_{\rm t}\equiv \{(x,z)\in \RR^{2},\ z=1\}, \\
\widehat{\phi_1} =\widehat{\psi} & \mbox{ on } \Gamma\equiv \{(x,z)\in \RR^{2},\ z=0\},
\end{array}
\right.
\end{eqnarray}
One can easily remark that~\eqref{Laplace1four} is an ordinary linear differential equation of order 2 whose solution has  the following form 
(since $\mu k^2 \geq 0$):
\[\widehat{\phi_1}(k)=A(k)\cosh (\sqrt{\mu}|k|z)+B(k)\sinh (\sqrt{\mu}|k|z).\]
Writing the boundary condition at  $z=0$ i.e. 
$\widehat{\phi_1}|_{z=0} =\widehat{\psi} $ and at $z=1$ i.e. $(\partial_z \widehat{\phi_1})|_{z=1}=0$ we obtain:
$$ A(k)=\widehat{\psi}(k) \quad \mbox{and} \quad B(k)=-\tanh(\sqrt{\mu}|k|)\widehat{\psi}(k) \ .$$
Thus we obtain:
$$ \widehat{\phi_1}(k)=\Big[\cosh (\sqrt{\mu}|k|z)-\tanh(\sqrt{\mu}|k|)\sinh (\sqrt{\mu}|k|z)\Big]\widehat{\psi}(k)  \ . $$
Therefore as:
\[ \widehat{G^{\mu}[0]\psi}=(\partial_z \widehat{\phi_1})|_{z=0} = -\sqrt{\mu}|k|\tanh(\sqrt{\mu}|k|)\widehat{\psi}(k),\]
we finally have:
\begin{equation}\label{Gmu0} G^{\mu}[0]\psi= -\sqrt{\mu}|k|\tanh(\sqrt{\mu}|k|)\psi(k).\end{equation}
$\bullet$ Linearizing~\eqref{Laplace2} around the rest state ($\zeta=0$) gives:
\begin{eqnarray}
\label{Laplace2lin} &\left\{
\begin{array}{ll}
\left(\ \mu\partial_x^2\ + \ \partial_z^2\ \right)\ \phi_2=0 & \mbox{ in } \Omega_2\equiv\{(x,z)\in \RR^{2},\ -\frac{1}{\delta}<z<0\}, \\
\partial_{z}\phi_2 = \partial_{z}\phi_1= G^{\mu}[0]\psi & \mbox{ on } \Gamma\equiv \{(x,z)\in \RR^{2},\ z=0\},\\
\partial_{z}\phi_2 =0 & \mbox{ on } \Gamma_{\rm b}\equiv \{(x,z)\in \RR^{2},\ z=-\frac{1}{\delta}\}.
\end{array}
\right.
\end{eqnarray}
Applying the Fourier transform with respect to $x$, one has:
\begin{eqnarray}
\label{Laplace2four} &\left\{
\begin{array}{ll}
-\mu|k|^2\widehat{ \phi_2}+\partial_z^2\widehat{ \phi_2}=0 & \mbox{ in } \Omega_2\equiv\{(x,z)\in \RR^{2},\ -\frac{1}{\delta}<z<0\}, \\
\partial_{z}\widehat{\phi_2} = \widehat{G^{\mu}[0]\psi}= -\sqrt{\mu}|k|\tanh(\sqrt{\mu}|k|)\widehat{\psi(k)}& \mbox{ on } \Gamma\equiv \{(x,z)\in \RR^{2},\ z=0\},\\
\partial_{z}\widehat{\phi_2} =0 & \mbox{ on } \Gamma_{\rm b}\equiv \{(x,z)\in \RR^{2},\ z=-\frac{1}{\delta}\}.
\end{array}
\right.
\end{eqnarray}
One can easily remark that~\eqref{Laplace2four} is an ordinary differential equation of order 2 whose solution has  the following form (since $\mu k^2 \geq 0$):
\[\widehat{\phi_2}(k)=A(k)\cosh (\sqrt{\mu}|k|z)+B(k)\sinh (\sqrt{\mu}|k|z).\]
Writing the boundary condition at  $z=0$ i.e. 
$(\partial_z \widehat{\phi_2})|_{z=0}=-\sqrt{\mu}|k|\tanh(\sqrt{\mu}|k|)\widehat{\psi(k)}$  and at $z=-1/\delta$ i.e. $(\partial_z \widehat{\phi_2})|_{z=-1/\delta}=0$ we obtain:
$$ B(k)=-\tanh(\sqrt{\mu}|k|)\widehat{\psi}(k) \quad \mbox{and} \quad 
-A(k)\sinh\Big(\dfrac{\sqrt{\mu}|k|}{\delta}\Big)+B(k)\cosh\Big(\dfrac{\sqrt{\mu}|k|}{\delta}\Big)=0 \ .$$
Replacing $B(k)$ by its expression, we obtain:
\[A(k)=-\dfrac{\tanh\Big(\sqrt{\mu}|k|\Big)}{\tanh\Big(\dfrac{\sqrt{\mu}|k|}{\delta}\Big)}\widehat{\psi}(k).\]
Since $\widehat{\phi_2}(k)|_{z=0}=A(k)$ we have thus:
 $\phi_2(k)|_{z=0}=-\dfrac{\tanh\Big(\sqrt{\mu}|k|\Big)}{\tanh\Big(\dfrac{\sqrt{\mu}|k|}{\delta}\Big)}\psi(k).$\\
Therefore,
\begin{equation}\label{Hmu0}H^{\mu,\delta}[0]\psi=\partial_x(\phi_2|_{z=0})=
-\dfrac{\tanh\Big(\sqrt{\mu}|k|\Big)}{\tanh\Big(\dfrac{\sqrt{\mu}|k|}{\delta}\Big)}\partial_x\psi .\end{equation}
Having obtained the two operators  $G^{\mu}$ and $H^{\mu,\delta}$ at the rest state, \eqref{Gmu0} and \eqref{Hmu0}, we can now proceed to compute the dispersion relation for the \emph {full Euler system}.\\ 
One derives the dispersion relation associated to~\eqref{Eulercomplet} by looking for plane wave solutions of the form  $(\zeta,\psi) = (\zetabar,\psibar)e^{i(kx-wt)}$ 
to the linearized equations around the rest state $(\zeta,\psi )=(0,0)$.\\
\\
The linearisation of the \emph {full Euler system}, \eqref{Eulercomplet} at the rest state, $\zeta = 0 \,,\, \psi = 0$,  writes under the form:
\begin{equation}\label{Eulercompletlin}
\left\{ \begin{array}{l}
\displaystyle\partial_{ t}{\zeta} \ -\ \frac{1}{\mu}G^{\mu}[0]\psi\ =\ 0,  \\ \\
\displaystyle\partial_{ t}\Big(H^{\mu,\delta}[0]\psi-\gamma \partial_x{\psi} \Big)\ + \ (\gamma+\delta)\partial_x{\zeta} = -\mu\frac{\gamma+\delta}{\bo}\frac{\partial_x \big(k(\epsilon\sqrt\mu\zeta)\big)}{{\epsilon\sqrt\mu}} \ ,
\end{array}
\right.
\end{equation}
Replacing $G^\mu[0]\psi$ by its expression given by \eqref{Gmu0} in the first equation of the system \eqref{Eulercompletlin}, we firstly have:
\[-iw\zetabar+\dfrac{1}{\sqrt{\mu}}|k|\tanh\Big(\sqrt{\mu}|k|\Big)\psibar=0 .\]
Thus, 
\begin{equation}\label{zetabar}\zetabar=\dfrac{1}{\sqrt{\mu}iw}|k|\tanh\Big(\sqrt{\mu}|k|\Big)\psibar.\end{equation}
Replacing $H^\mu[0]\psi$ by its expression given by \eqref{Hmu0} in the second equation of the system~\eqref{Eulercompletlin}, we obtain:
\[-\dfrac{\tanh\Big(\sqrt{\mu}|k|\Big)}{\tanh\Big(\dfrac{\sqrt{\mu}|k|}{\delta}\Big)}\partial_t(\partial_x\psi)-\gamma\partial_t(\partial_x\psi)+(\gamma+\delta)\partial_x\zeta-\mu\frac{\gamma+\delta}{\bo}\partial_x^3\zeta=0 .\]
Thus,
\[-\dfrac{\tanh\Big(\sqrt{\mu}|k|\Big)}{\tanh\Big(\dfrac{\sqrt{\mu}|k|}{\delta}\Big)}wk\psibar-
\gamma wk \psibar +(\gamma+\delta)ik\zetabar+\mu\frac{\gamma+\delta}{\bo}ik^3\zetabar=0.\]
After straightforward computations and replacing $\zetabar$ by its expression given in~\eqref{zetabar} we finally obtain the following dispersion relation:
\begin{equation}\label{rdfe}
w^2_{F.E}=\dfrac{(\gamma+\delta)|k|\tanh\Big(\sqrt{\mu}|k|\Big)(1+\dfrac{\mu}{\bo}k^2)\tanh\Big(\sqrt{\mu}\dfrac{|k|}{\delta}\Big)}
{\sqrt{\mu}[\tanh\Big(\sqrt{\mu}|k|\Big)+\gamma\tanh\Big(\sqrt{\mu}\dfrac{|k|}{\delta}\Big)]}.\end{equation}
%
%
\subsubsection{Phase velocity agreement}
\begin{figure}[H]
\centering
\includegraphics[scale=1]{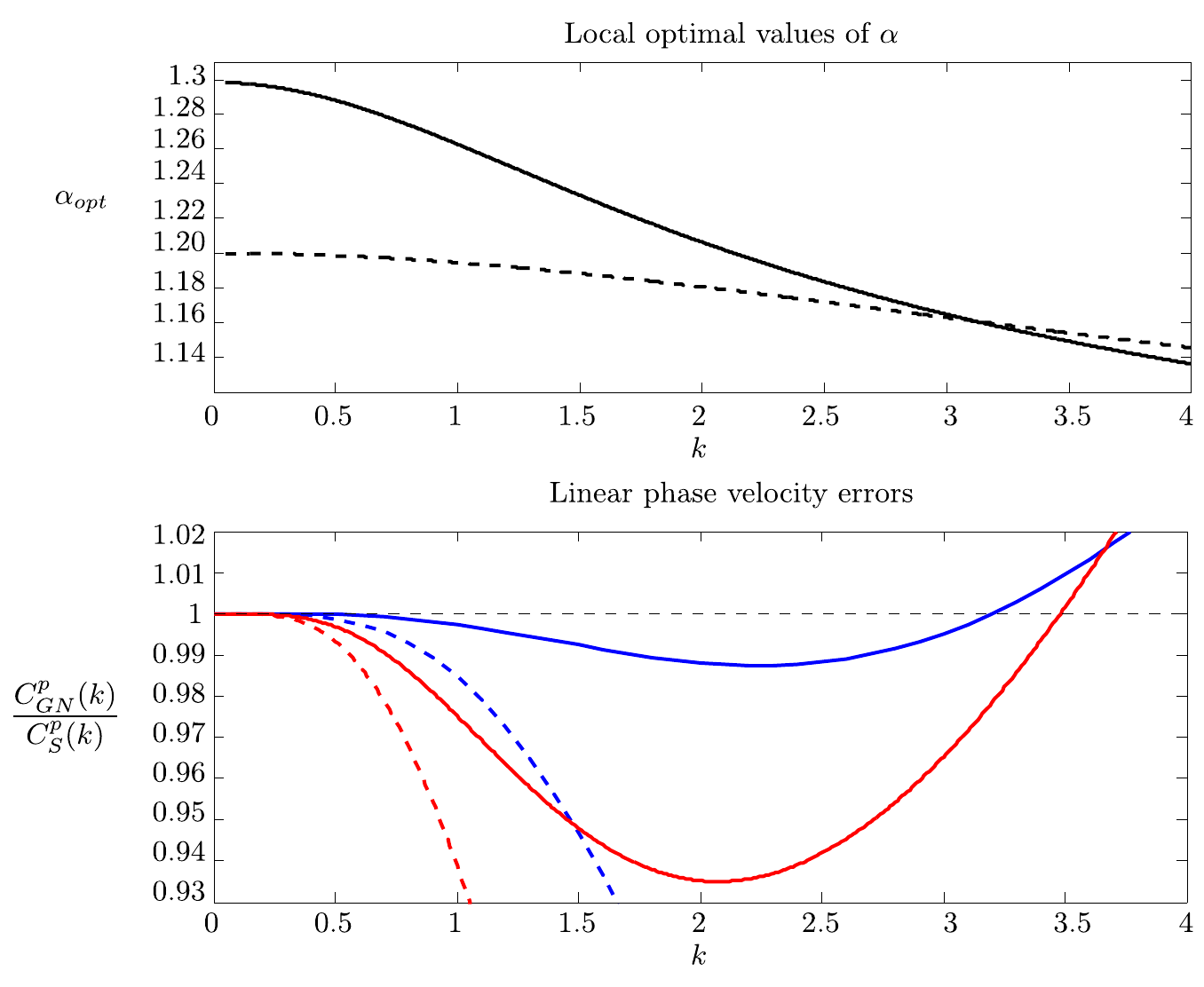}
\caption{Top: local optimal values of $\alpha$ for $k\in[0,4]$, for the two layers case (full line) and one layer case (dashed line). Bottom: linear phase velocity errors in the one layer case (blue) and two layers case (red) for $\alpha=1.159$ (full blue line) and $\alpha=1.498$ (full red line) and $\alpha=1$ (dashed line) }
\label{figalphaopt}
\end{figure}

One can easily remark that the Taylor expansions of the two previous dispersion relations~\eqref{rdgn} and~\eqref{rdfe} are equivalent for small wavenumbers,
\[w_{\alpha,GN}^2\equiv w^2_{F.E}\simeq k^2-\mu\nu k^4 + \OO(\mu^2 k^6).\]
Indeed, for small wavenumber one has:
\begin{eqnarray*}w^2_{\alpha,GN}&=&k^2\Big(\dfrac{1+\mu\nu(\alpha-1)k^2}{1+\mu\nu\alpha k^2}\Big)
\\&\approx& k^2(1+\mu\nu(\alpha-1)k^2)(1-\mu\nu\alpha k^2+\OO((\sqrt{\mu}k)^4))
\\&\approx &k^2(1-\mu\nu\alpha k^2+\mu\nu(\alpha-1)k^2+\OO((\sqrt{\mu}k)^4))
\\&\approx&k^2-\mu\nu k^4+\OO(\mu^2 k^6).\end{eqnarray*} 
and
\begin{eqnarray*}
w^2_{FE}&=&\dfrac{(\gamma+\delta)|k|\tanh(\sqrt{\mu}|k|)(1+\dfrac{\mu}{\bo}k^2)\tanh(\sqrt{\mu}\dfrac{|k|}{\delta})}{\sqrt{\mu}[\tanh(\sqrt{\mu}|k|)+\gamma\tanh(\sqrt{\mu}\dfrac{|k|}{\delta})]}
\\&\approx& \dfrac{(\gamma+\delta)|k|\big[\sqrt{\mu}k-\dfrac{(\sqrt{\mu}k)^3}{3}+\OO((\sqrt{\mu}k)^5)\big](1+\dfrac{\mu}{\bo}k^2)\big[\dfrac{\sqrt{\mu}k}{\delta}-\big(\dfrac{\sqrt{\mu}k}{\delta}\big)^3\dfrac{1}{3}+\OO((\sqrt{\mu}k)^5)\big]}{\sqrt{\mu}\Big[\sqrt{\mu}k-\dfrac{(\sqrt{\mu}k)^3}{3}+\OO((\sqrt{\mu}k)^5)+\gamma\Big(\dfrac{\sqrt{\mu}k}{\delta}-\big(\dfrac{\sqrt{\mu}k}{\delta}\big)^3\dfrac{1}{3}+\OO((\sqrt{\mu}k)^5)\Big)\Big]}
\\&\approx&k^2+\mu k^4\Big[\dfrac{1}{\bo}-\dfrac{\delta^2(1+\gamma\delta)}{3\delta^3(\gamma+\delta)}\Big]+\OO(\mu^2 k^6)
\\&\approx&k^2-\mu\nu k^4+\OO(\mu^2 k^6).
\end{eqnarray*}
Therefore for small wavenumbers and in the whole range of regime, the choice of $\alpha$ does not influence the dispersion relation \eqref{rdgn}.
In the desired simulations, we are interested in small wave length (i.e large wavenumbers) dispersion characteristics and thus we would like to find an optimal value of $\alpha$ in order
to observe the same dispersion properties for the reduced Green-Naghdi model that the \emph {full Euler} original problem satisfies. Therefore, we are interested in values of $\alpha$
such that $C_{GN}^p(k)=C_{F.E}^p(k)$ for large value of $k$.
We can thus compute $\alpha_{opt}$ from the equality: $w_{\alpha,GN}^2= w^2_{F.E}$.\\
Let us denote $X=\sqrt{\mu}|k|$. The previous equality writes:
\[\dfrac{1+\nu X^2(\alpha_{opt}-1)}{1+\nu X^2\alpha_{opt}}
=\dfrac{(\gamma+\delta)\tanh(X)(1+\dfrac{X^2}{\bo})\tanh(\dfrac{X}{\delta})}{X[\tanh(X)+\gamma\tanh(\dfrac{X}{\delta})]}=g(X) .\]
After straightforward computations we obtain the follwing expression for the value of $\alpha_{opt}$:
\[\alpha_{opt}=\dfrac{g(X)-1+\nu X^2}{\nu X^2(1-g(X))}.\]
In Figure~\ref{figalphaopt} (top), $\alpha_{opt}$ is plotted against the spatial wave number $k$, for $k \in [0,4]$ and $\mu=1$, for the one and two layers cases. For the one layer case we set $\gamma=0$, $\delta=1$ and $\bo^{-1}=0$ whereas for the two layers case we set $\gamma=0.95$, $\delta=0.5$ and $\bo^{-1}=5\times 10^{-5}$. We would like to mention that $\alpha_{opt}\rightarrow 1$ in both cases when considering very large values of $k$, this fact is confirmed by the numerical simulation done in Section~\ref{DAM2}, see Figure~\ref{dam2}, where we notice similar observations at time $t=20 \ s$.  Therefore, we will use an optimal value of $\alpha$ different of $1$ only when looking at the dispersion effects for intermediate wave numbers. This is highlighted, when comparing to the numerical experiments done in~\cite{DucheneIsrawiTalhouk16} (for more details see Section~\ref{KH}).
In Figure~\ref{figalphaopt} (bottom), the ratio $\dfrac{C^p_{GN}(k)}{C^p_{S}(k)}$ (i.e linear phase velocity error), is plotted against the spatial wave number $k$ for $k\in[0,4]$, in the one and two layer cases.

\subsection{High frequencies instabilities}\label{IHFsec}
\begin{figure}[H]
    \centering
{\includegraphics[scale=1]{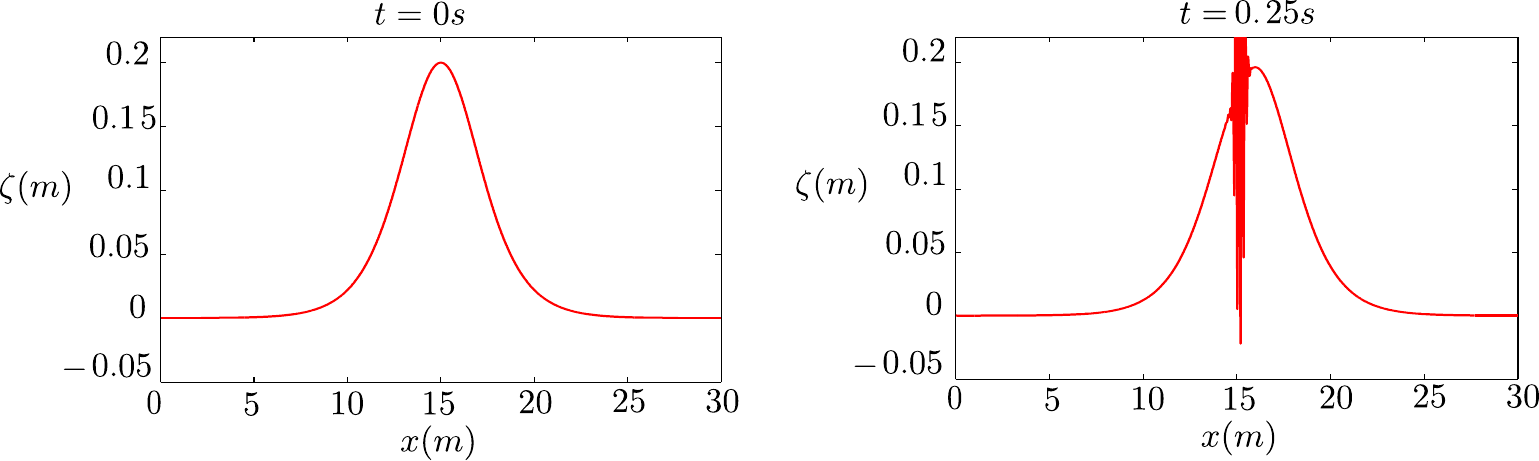}}
    \caption{High frequencies instabilities in the one layer case due to third order derivatives.}
    \label{figIHF}
\end{figure}

In what follows, a qualitative discussion on the stability of the improved model~\eqref{GNCH6} is given. We study the high frequencies instabilities in the one layer case, as well as in the two layers case.
As $\mu <<1$, a simple expansion of the operator  $(I+\mu\nu \alpha T[0])^{-1}$ shows that we have: 
$$(I+\mu\nu \alpha T[0])^{-1}[(\gamma+\delta)\partial_x \zeta ]=(\gamma+\delta)\partial_x \zeta +\OO(\mu).$$ Thus, we could replace 
$ Q_2(\zeta)$ by $\widetilde{Q}_2(\zeta)=-S[\zeta]\big((\gamma+\delta)\partial_x \zeta\big),$ and 
$ Q_3(\zeta)$ by $\widetilde{Q}_3(\zeta) =  \kappa_1\zeta T[0][(\gamma+\delta)\partial_x \zeta ]$
in the second equation of \eqref{GNCH6}, keeping the same order of precision $\OO(\mu^2)$.
This simple expansion would avoid the inversion of $(I+\mu\nu \alpha T[0])$ (resolution of an extra linear system) in the computation of
${Q}_2$ and ${Q}_3$ but it leads to instabilities.
Indeed, this is due to the fact that the two terms  $\widetilde{Q}_2(\zeta)$ and  $\widetilde{Q}_3(\zeta)$
contain third order derivatives in $\zeta$ that may create high frequencies instabilities.
\begin{remark}
In order to recover the dimensionalized version of system~\eqref{GNCH6}, one has to set $\mu=\epsilon=1$ and add the acceleration of gravity term $g$ when necessary to obtain the following system:
\begin{equation}\label{GNCH6dim}\left\{ \begin{array}{l}
        \dsp \partial_{ t}\zeta +\partial_x\big(f(\zeta)  v\big)\ =\  0,\\ \\
        \dsp (I+\nu\alpha T[0]) \big[ \partial_{ t}  v+\varsigma v\partial_x v+\dfrac{\alpha-1}{\alpha}\big((\gamma+\delta)g \partial_x \zeta +\partial_x(q_3(\zeta)  {v}^2)\big)\big]\\+\dfrac{1}{\alpha}\big((\gamma+\delta)g\partial_x \zeta +\partial_x(q_3(\zeta)  {v}^2)\big)+Q_1(v)+\nu Q_2(\zeta)+\nu Q_3(\zeta)=0,
    \end{array} \right. \end{equation}
with
\begin{equation}\label{Q2dim2layers}Q_2(\zeta)= \kappa_2 \partial_x \Big(\zeta \partial_x \big((I+\nu\alpha T[0])^{-1}[(\gamma+\delta)g\partial_x \zeta ]\big)\Big)\end{equation}
and
\begin{equation}\label{Q3dim2layers}Q_3(\zeta)=\kappa_1\zeta T[0](I+\nu\alpha T[0])^{-1}[(\gamma+\delta)g\partial_x \zeta ].\end{equation}
\end{remark}
Firstly, we discuss the stability of the model~\eqref{GNCH6dim} in the one layer case without surface tension. To this end, we set $\gamma=0$, $\delta=1$ and $\dfrac{1}{\bo}=0$. Thus system~\eqref{GNCH6dim} becomes:
\begin{equation}\label{GNCH6dim1layer}\left\{ \begin{array}{l}
        \dsp \partial_{ t}\zeta +\partial_x\big((1+\zeta)  v\big)\ =\  0,\\ \\
        \dsp (I+\frac{\alpha}{3} T[0]) \big[ \partial_{ t}  v+v\partial_x v+\dfrac{\alpha-1}{\alpha}g \partial_x \zeta\big]+\dfrac{1}{\alpha}g\partial_x \zeta+Q_1(v)+ Q_2(\zeta)+\dfrac{1}{3} Q_3(\zeta)=0,
    \end{array} \right. \end{equation}
with
\begin{equation*}
    Q_1(v)=\dfrac{2}{3}\partial_x ((\partial_x v)^2),
\end{equation*}
\begin{equation}\label{Q2dim1layer}   
    Q_2(\zeta)= \partial_x \Big(\zeta \partial_x\big((I+\frac{\alpha}{3} T[0])^{-1}[g\partial_x \zeta]\big)\Big),
\end{equation}
    \begin{equation}\label{Q3dim1layer}
    Q_3(\zeta)=\zeta T[0](I+\frac{\alpha}{3}T[0])^{-1}[g\partial_x \zeta ].
\end{equation}
When linearizing system~\eqref{GNCH6dim1layer} around constant state solution $(\underline{\zeta}, \underline{v})$, one obtains the following linear system in $(\tilde{\zeta},\tilde{v})$:
\begin{equation}\label{GNCH6dim1layerlin}\left\{ \begin{array}{l}
        \dsp \partial_{ t}\tilde{\zeta} +(1+\underline{\zeta}) \partial_x \tilde{v} + \underline{v} \partial_x{\tilde{\zeta}}\ =\  0,\\ \\
        \dsp (I-\frac{\alpha}{3} \partial_x^2) \big[ \partial_{ t}  \tilde{v}+\underline{v}\partial_x \tilde{v}+\dfrac{\alpha-1}{\alpha}g \partial_x \tilde{\zeta}\big]+\dfrac{1}{\alpha}g\partial_x \tilde{\zeta}+\dfrac{2}{3}g\underline{\zeta}\partial_x^2\Big((I-\frac{\alpha}{3} \partial_x^2)^{-1}\partial_x \tilde{\zeta}\Big)=0.
    \end{array} \right. \end{equation}
Looking for plane wave solution of the form $(\tilde{\zeta},\tilde{v})=e^{i(kx-wt)}(\zeta^0,v^0)$ as solution of the above system, one obtains the following dispersion relation:
\begin{equation}\label{rd}
    \dfrac{(w-k\underline{v})^2}{g(1+\underline{\zeta}) k^2}=\dfrac{1+\frac{(\alpha-1)}{3}k^2-\dfrac{2k^2\underline{\zeta}}{3(1+\frac{\alpha}{3} k^2)}}{1+\frac{\alpha}{3} k^2}.
\end{equation}
When we consider the linearization of our new model around the rest state $(\underline{\zeta},\underline{v})=(0, 0)$, the dispersion relation~\eqref{rd} becomes:
\begin{equation*}\label{rd0}
    w^2=gk^2\dfrac{1+\frac{(\alpha-1)}{3}k^2}{1+\frac{\alpha}{3} k^2}.
\end{equation*}
In this case the perturbations are always stable if $\alpha \geq 1$.
We refer to Section~\ref{Secalphachoice} for the discussion concerning the choice of $\alpha$ in order to improve the dispersive properties of the model.\\
Now, replacing $ Q_2(\zeta)$ defined in~\eqref{Q2dim1layer} by $\widetilde{Q}_2(\zeta)=\partial_x \Big(\zeta \partial_x\big(g\partial_x \zeta\big)\Big)$ and 
$  Q_3(\zeta)$ defined in~\eqref{Q3dim1layer} by $\widetilde{Q}_3(\zeta) = \kappa_1\zeta T[0][g\partial_x \zeta ]$ in the second equation of \eqref{GNCH6dim1layer}, 
modifies the dispersion relation~\eqref{rd} and we obtain instead:
\begin{equation}\label{rdIHF}
    \dfrac{(\tilde{w}-k\underline{v})^2}{g(1+\underline{\zeta}) k^2}=\dfrac{1+\frac{(\alpha-1)}{3}k^2-\frac{2k^2\underline{\zeta}}{3}}{1+\frac{\alpha}{3} k^2}.
\end{equation}
One can easily check that if $\underline{\zeta} > \frac{(\alpha -1)}{2}$, the numerator of the right hand side of  the dispersion relation~\eqref{rdIHF} will become negative for sufficiently large values of $k^2$. Thus the root of $\tilde{w}$ will become complex inducing a high frequency instability of the model, see Figure~\ref{figIHF}. On the other hand, the numerator of the r.h.s of the dispersion relation~\eqref{rd} is positive for sufficiently large values of $k^2$, due to the existence of the term $-\dfrac{2 k^2 \underline{\zeta}}{3(1+\frac{\alpha}{3}k^2)}$, thus $w$ is always real at high frequencies.

\begin{figure}[H]
    \centering
{\includegraphics[scale=0.99]{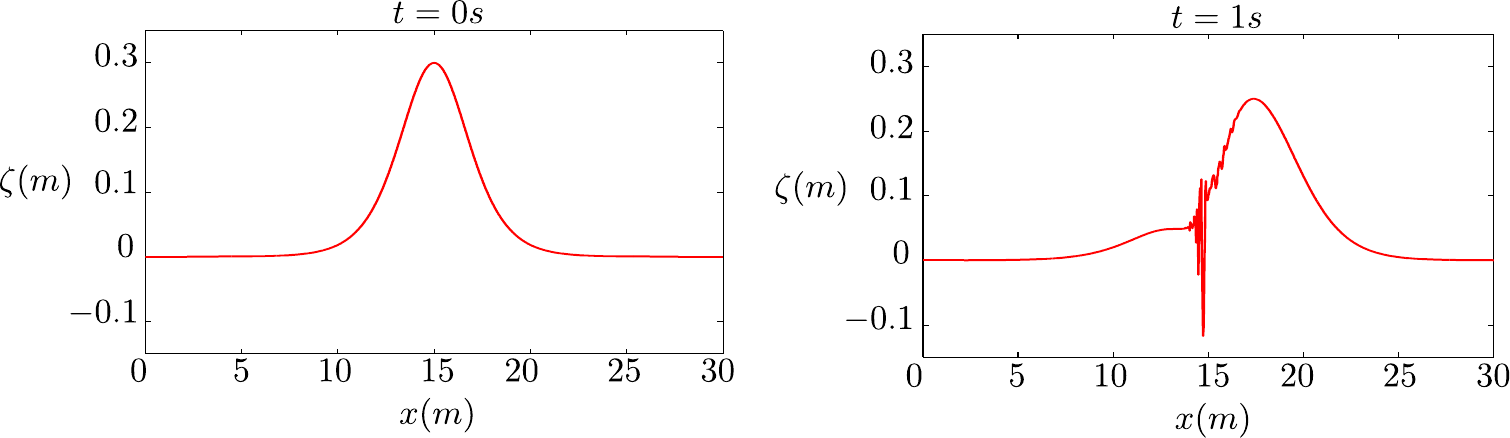}}
    \caption{High frequencies instabilities in the two layers case due to third order derivatives.}
    \label{figIHF2layers}
\end{figure}

Let us now discuss the stability issue of the model~\eqref{GNCH6dim} in the two layers case without surface tension. When linearizing system~\eqref{GNCH6dim} around motionless steady state solution $(\underline{\zeta}=\text{cst}, \underline{v}=0)$ and after following the same method as above, one obtains the following dispersion relation:
\begin{equation}\label{rd2layers}
    \dfrac{w^2}{gf(\underline{\zeta}) k^2}=\dfrac{(\gamma+\delta)\Big(1+\nu(\alpha-1)k^2-\nu\dfrac{(\kappa_2-\kappa_1)k^2\underline{\zeta}}{(1+\nu\alpha k^2)}\Big)}{1+\nu\alpha k^2}.
\end{equation}
Replacing  $Q_2(\zeta)$ defined in~\eqref{Q2dim2layers} by $\widetilde{Q}_2(\zeta)=\partial_x \Big(\zeta \partial_x\big(g(\gamma+\delta)\partial_x \zeta\big)\Big)$ and 
$Q_3(\zeta)$ defined in~\eqref{Q3dim2layers} by $\widetilde{Q}_3(\zeta) = \kappa_1\zeta T[0][g(\gamma+\delta)\partial_x \zeta ]$ in the second equation of \eqref{GNCH6dim}, 
modifies the dispersion relation~\eqref{rd2layers} and we obtain instead:
\begin{equation}\label{rd2layersIHF}
    \dfrac{\tilde{w}^2}{gf(\underline{\zeta}) k^2}=\dfrac{(\gamma+\delta)(1+\nu(\alpha-1)k^2-\nu(\kappa_2-\kappa_1)k^2\underline{\zeta})}{(1+\nu\alpha k^2)}.
\end{equation}
In this case, there exists a critical ratio for the depth of the two layers. Indeed, when $\delta^2 < \gamma$, one has $\kappa_2<\kappa_1$, thus the perturbations are always stable if $\alpha \geq 1$. Whereas, when $\delta^2\geq\gamma$ i.e assuming $\gamma=0.5$, $\delta=0.8$ and $\dfrac{1}{\bo}=0$, one has $\kappa_2>\kappa_1$, thus the perturbations are unstable if $\underline{\zeta} > \dfrac{\alpha -1}{\kappa_2-\kappa_1}$. In fact, under the previous conditions the numerator of the right hand side of  the dispersion relation~\eqref{rd2layersIHF} will become negative for sufficiently large values of $k^2$. Thus the root of $\tilde{w}$ will become complex inducing a high frequency instability of the model, see Figure~\ref{figIHF2layers}. On the other hand, the numerator of the r.h.s of the dispersion relation~\eqref{rd2layers} is positive for sufficiently large values of $k^2$, due to the existence of the term $-\nu\dfrac{(\kappa_2-\kappa_1)k^2\underline{\zeta}}{(1+\nu\alpha k^2)}$, thus $w$ is always real at high frequencies.  This ensures the numerical stability of the model~\eqref{GNCH6dim} which will be considered in the rest of this work.
\section{Numerical methods}\label{NMSec}
This section is devoted to the numerical methods developped to solve the improved Green-Naghdi equations~\eqref{GNCH6}.
As pointed out by many authors~\cite{BCLMT,LannesMarche14} the improved dispersion Green-Naghdi equations~\eqref{GNCH6} is well-adapted to the implementation of a splitting scheme separating the hyperbolic and the dispersive parts of the equations. 
We present in Section~\ref{SSsec} this splitting scheme inspired by~\cite{BCLMT,LannesMarche14}.
We explain in Sections~\ref{Hyperbolic} and~\ref{Dispersive} how we treat respectively the hyperbolic and dispersive parts of the equations.

We decided to treat the hyperbolic part by a finite volume method of Roe type. We will construct first order, second order ``MUSCL'' type method and finally 5th order WENO method. The high order method is suitable to compute correctly the maximum value of the height $\zeta$ and the discontinuities
by limiting the diffusive effects.
As we will show in the numerical validations, the high order scheme is also suitable to catch correctly the dispersive effects.
 
The dispersive part of the proposed splitting method is solved using a classical finite difference method.

\subsection{The splitting scheme}\label{SSsec}
Let us recall the improved Green-Naghdi system that we consider:
\begin{equation}\label{GNCH66}\left\{ \begin{array}{l}
        \dsp \partial_{ t}\zeta +\partial_x\big(f(\epsilon\zeta)  v\big)\ =\  0,\\ \\
        \dsp (I+\mu\nu\alpha T[0]) \big[ \partial_{ t}  v+\epsilon\varsigma v\partial_x v+\dfrac{\alpha-1}{\alpha}\big((\gamma+\delta)\partial_x \zeta +\epsilon\partial_x(q_3(\epsilon\zeta)  {v}^2)\big)\big]\\+\dfrac{1}{\alpha}\big((\gamma+\delta)\partial_x \zeta +\epsilon\partial_x(q_3(\epsilon\zeta)  {v}^2)\big)+\mu\epsilon Q_1(v)+\mu\nu Q_2(\zeta)+\mu\epsilon\nu Q_3(\zeta) =0.
    \end{array} \right. \end{equation}
$q_3$ defined by~\eqref{defkappaq3},$Q_1$, $Q_2$ and $Q_3$ are defined by~\eqref{Q1},~\eqref{Q2} and~\eqref{Q3}.\\

We decompose the solution operator $S(.)$ associated to~\eqref{GNCH66} at each time step $\Delta t$ by the following second order splitting scheme:
$$S(\Delta t) = S_1(\Delta t/2)S_2(\Delta t)S_1(\Delta t/2)$$
where $S_1(.)$ is the solution operator associated to the hyperbolic part, and $S_2(.)$ the solution operator associated to the dispersive part of the
equations \eqref{GNCH66}.

$\bullet \ S_1(t)$ is the solution operator associated to the hyperbolic part namely the nonlinear shallow water equations, NSWE:
\begin{equation}\label{hyp}\left\{ \begin{array}{l}
        \dsp \partial_{ t}\zeta +\partial_x\big(f(\epsilon\zeta)  v\big)\ =\  0,\\ \\
        \dsp \partial_{ t}  v+\epsilon\varsigma v\partial_x v+\dfrac{\alpha-1}{\alpha}\big((\gamma+\delta)\partial_x \zeta +\epsilon\partial_x(q_3(\epsilon\zeta)  {v}^2)\big)+\dfrac{1}{\alpha}\big((\gamma+\delta)\partial_x \zeta +\epsilon\partial_x(q_3(\epsilon\zeta)  {v}^2)\big) =0.
    \end{array} \right. \end{equation}
Using the definition of $q_3(\epsilon\zeta)$ given in~\eqref{defkappaq3}, one can easily check that 
$\dsp{\dfrac{\varsigma}{2}+q_3(\epsilon\zeta)=\dfrac {f'(\epsilon\zeta)}{2}}$. Thus we rewrite the NSWE system \eqref{hyp} in the following condensed form:
\begin{equation}\label{hypcons}\left\{ \begin{array}{l}
        \dsp \partial_{ t}\zeta +\partial_x\big(f(\epsilon\zeta)  v\big)\ =\  0,\\ \\
        \dsp \partial_{ t}  v+\partial_x\Big(\dfrac{\epsilon f'(\epsilon\zeta)}{2}v^2+(\gamma+\delta)\zeta\Big) =0.
    \end{array} \right. \end{equation}
We recall that, $f(\epsilon\zeta) =\dsp\frac{h_1h_2}{h_1+\gamma h_2}$ and $f'(\epsilon\zeta)  = \dsp\frac{h_1^2-\gamma h_2^2}{(h_1+\gamma h_2)^2}$, with $h_1=1-\epsilon\zeta$ and $h_2=1/\delta+\epsilon\zeta$.

$\bullet \ S_2(t)$ is the solution operator associated to the remaining (dispersive) part of the equations.
\begin{equation}\label{disp}\left\{ \begin{array}{l}
        \dsp \partial_{ t}\zeta \ =\  0,\\ \\
        \dsp (I+\mu\nu\alpha T[0]) \big[ \partial_{ t}  v-\dfrac{1}{\alpha}\big((\gamma+\delta)\partial_x \zeta +\epsilon\partial_x(q_3(\epsilon\zeta)  {v}^2)\big)\big]\\+\dfrac{1}{\alpha}\big((\gamma+\delta)\partial_x \zeta +\epsilon\partial_x(q_3(\epsilon\zeta)  {v}^2)\big)+\mu\epsilon Q_1(v)+\mu\nu Q_2(\zeta)+\mu\epsilon\nu Q_3(\zeta)=0.
    \end{array} \right. \end{equation}
 In this study, $S_1$ is computed using a finite volume method while $S_2$ is computed using a classical finite-difference  method.
 
\begin{figure}[H]
     \begin{center}
         \includegraphics[height = 2cm]{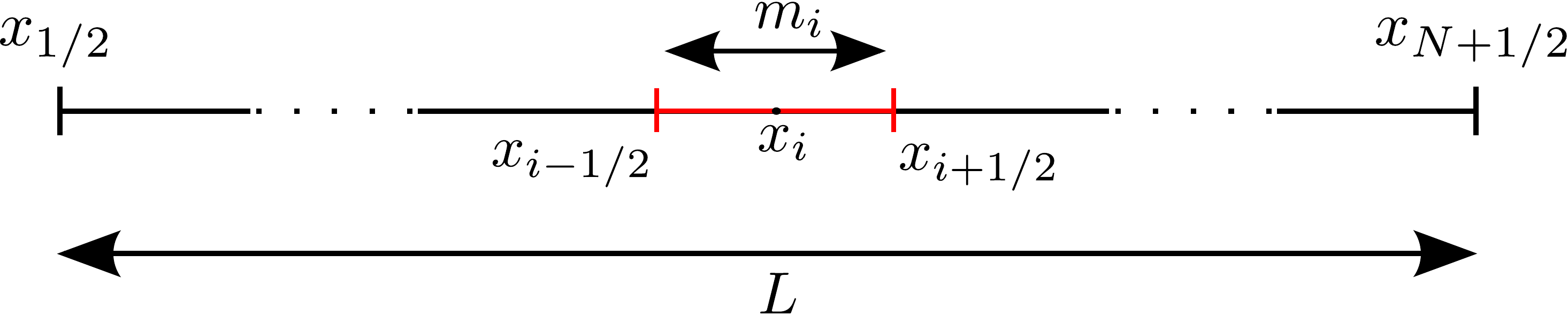}
         \caption{The space discretization.}
         \label{domain-discret}
        \end{center}
    \end{figure}

 In order to discretize system \eqref{GNCH66}, the numerical domain is an interval of length $L$ denoted $[0,L]$. 
Let $N\in\NN^{*}$, and let us consider the following mesh on $[0,L]$. Cells are denoted for every
 $i\in [0,N+1]$, by $m_i =(x_{i-1/2},x_{i+1/2})$, with  $x_i=\displaystyle\frac{x_{i-1/2}+x_{i+1/2}}{2}$ and 
 $h_{i}=x_{i+1/2}-x_{i-1/2} $ the space mesh. The ``fictitious'' cells $m_0$ and $m_{N+1}$ denote the boundary cells and 
 the mesh  interfaces located at $x_{1/2} = 0$ and $x_{N+1/2} = L$ are respectively the upstream and the downstream ends  (see Figure \ref{domain-discret}).
 
 We denote $\displaystyle \Delta x = \min_{i=1,N}h_i$.
 
    We also consider a time discretization $t^n$ defined by $t^{n+1}=t^n+\Delta t$ with $\Delta t$ the time step. 
    
    
\subsection{Finite volume scheme}\label{Hyperbolic}
For the sake of simplicity in the notations, it is convenient to rewrite the hyperbolic system~\eqref{hypcons} in the following form:
\begin{equation}\label{condform}\partial_t U +\partial_x (F(U))=0,\end{equation}
with the following conservative variables and flux function:
\begin{equation}\label{consvar}
U=\begin{pmatrix}
\zeta\\
v
\end{pmatrix}
,\quad F(U)=\begin{pmatrix}
f(\epsilon\zeta)v\\
\dfrac{\epsilon f'(\epsilon\zeta)}{2} v^2+(\gamma+\delta)\zeta
\end{pmatrix}.
\end{equation}
The Jacobian matrix is given by:
\begin{equation}\label{matricejacob}
A(U)=d(F(U))=\begin{pmatrix}
\epsilon f'(\epsilon\zeta)v&f(\epsilon\zeta)\\
(\gamma+\delta)+\epsilon^2\frac{f''(\epsilon\zeta)}{2}v^2 & \epsilon f'(\epsilon\zeta)v
\end{pmatrix},
\end{equation}
where $f''(\epsilon\zeta)= -\dfrac{2\gamma(h_1+h_2)^2}{(h_1+\gamma h_2
)^3}$.\\
\\
A simple computation shows that the homogeneous system~\eqref{condform} is strictly hyperbolic provided that:
\begin{equation}\label{Condhyp}\left\{ \begin{array}{l}
\inf_{x\in \RR} h_1 >\ 0,\\ \\
\inf_{x\in \RR} h_2 >\ 0,\\ \\
\inf_{x\in \RR} \Big[(\gamma+\delta) -\dfrac{\gamma(h_1+h_2)^2}{(h_1+\gamma h_2
)^3}\epsilon^2v^2 \Big]>\ 0.\\ \\
\end{array} \right. \end{equation}
As a matter of fact, these conditions simply consist in assuming that the deformation of the interface is not too large and imposing a smallness assumption on $\epsilon v$. Notice that theses conditions correspond exactly to the hyperbolicity condition for the shallow water system provided in~\cite{GuyenneLannesSaut10}.\\

As a consequence, the solutions may develop shock discontinuities. In order to rule out the
unphysical solutions, the system~\eqref{condform} must be supplemented by entropy inequalities, (see for instance~\cite{Bouchut04} and references therein for more details).
The Cauchy problem associated to~\eqref{condform} is the following:
\begin{equation}\label{cauchy}\left\{ \begin{array}{l}
\partial_t U +\partial_x (F(U)) \ =\  0, \qquad  t\geq0, x\in \RR.\\ \\
U(0,x)=U_0(x), \qquad  x \in \RR.
\end{array} \right. \end{equation}
We are interested in the approximation of~\eqref{cauchy} by the finite volume method.\\

We denote 
$\overline{U}_i=(\zeta_i,v_i)$, the cell-centered approximation of $U$ on the cell $m_i$ at time $t$ given by:
$$\overline U_i = \frac{1}{h_i}\int_{m_i}U(t,x)\,dx \ .$$
The piecewise constant representation of $U$  is given by,   $U(t,x) = \displaystyle \overline U_i \mathds{1}_{m_i}(x)$.

We denote $\overline{U}_i^n=(\zeta_i^n,v_i^n)$, the cell-centered approximation of $U$ on the cell $m_i$ at time $t^n$ given by:
$$\overline U^n_i = \frac{1}{h_i}\int_{m_i}U(t^n,x)\,dx \ .$$

The spatial discretization of the homogeneous system \eqref{hyp} can be recast under the following classical semi-discrete finite-volume formalism:
\begin{equation}\label{volfinidiscret}
\frac{d \overline{U}_i(t)}{dt} + \frac{1}{h_i}\Big(\widetilde{F}(\overline U_i,\overline U_{i+1}) - \widetilde{F}(\overline U_{i-1},\overline U_{i}) \Big) = 0
\end{equation}
where $\widetilde{F}$ is a numerical flux function based on a conservative flux consistent with the homogeneous NSWE:
\begin{equation}\label{flux} F_{i+1/2}=\widetilde{F}(\overline U_i,\overline U_{i+1})\approx\dfrac{1}{h_i} \int_{m_i} F(U(t,x_{i+1/2}))dx.\end{equation}
\paragraph{VFRoe method.}
In what follows, we consider the numerical approximation of the hyperbolic system of conservation laws in the form of~\eqref{condform}. To this end, we adopt the VFRoe method (see~\cite{BGH00,GHN02,GHN03}) which is an approximate Godunov scheme. It relies on the exact 
resolution of the following linearized Riemann problem:
\begin{equation}\label{PRL}\left\{ \begin{array}{l}
\dsp \partial_{ t}U + \widetilde{A}(\overline U_i^n, \overline U_{i+1}^n)\partial_x U\ =\  0,\\ \\
\dsp U(0,x)=\left \{ \begin{array}{l}
\dsp \overline U^n_i \quad if \quad x<x_{i+1/2},\\ \\
\dsp \overline U_{i+1}^n \quad if \quad x>x_{i+1/2},
\end{array}\right.\end{array} \right. 
\end{equation}
where $\widetilde{A}(\overline U_i^n,\overline U_{i+1}^n)=A\left(\dfrac{\overline U_i^n+\overline U_{i+1}^n}{2}\right)$.\\
By solving the linearized Riemann problem we obtain $\overline U_{i+1/2}^*= U(x=x_{i+1/2},t=t_{n})$, the  interface value between two neighbouring cells.\\
In what follows, we will detail the choice of the numerical flux for different order of approximation. The only change is in the computation of the interface values $\overline U_{i+1/2}^*$ which depends on the right and the left states of the linearized Riemann problem. For the sake of simplicity, we will
still denote by $\widetilde{F}$ these numerical fluxes.
\paragraph{CFL condition}\label{CFLpar}
It is always necessary to impose what is called a CFL condition (for
Courant, Friedrichs, Levy) on the timestep to prevent the blow up of
the numerical values. It comes usually under the form
\begin{equation}\label{CFL}a_{i+1/2}\Delta t \leq \Delta x, \quad i=1,\ldots,N,\end{equation}
where $a_{i+1/2}=\dsp \max_{i\in[1, N]}(j=1,2, |\lambda_j(\widetilde U_i)|)$ and $\lambda_j(\widetilde U_i)$ are the eigenvalues of $A\big(\widetilde U_i=\dfrac{\overline U_i^n+ \overline U_{i+1}^n}{2}\big)$.\\
The restriction~\eqref{CFL} enables in practice to compute the timestep
at each time level $t_n$, in order to determine the new time level $t_{n+1} =
t_n + \Delta t$ (within this view, $\Delta t$ is not constant, it is computed in an
adaptive fashion).
\paragraph{Consistency.}
The numerical flux $\widetilde{F}(U_l,U_r)$ is called consistent with~\eqref{condform} if
\begin{equation}\label{consis}\widetilde{F}(U,U)=F(U) \quad \mbox{for all U}.\end{equation}



\subsubsection{First order finite-volume scheme}\label{FV1sec}
The semi-discrete equation \eqref{volfinidiscret} is discretised by an explicit Euler (in time) method to obtain :
\begin{equation}\label{volfini}
\overline U^{n+1}_i= \overline U^n_i-\dfrac{\Delta t}{h_i}(F_{i+1/2}^n-F_{i-1/2}^n),
\end{equation}
where the numerical flux is defined directly as the value of the exact flux at the interface value, namely: 
\begin{eqnarray}\label{numflux}
F_{i+1/2}^n=\widetilde{F}(\overline U_i^n,\overline U_{i+1}^n)=F(\overline U_{i+1/2}^*)\nn \\ 
\\F_{i-1/2}^n=\widetilde{F}(\overline U_{i-1}^n,\overline U_i^n)=F(\overline U_{i-1/2}^*)\nn.
\end{eqnarray}

Let us remark that by construction the numerical flux given by \eqref{numflux} ensures the consistency property.

In the sequel, we will suppose that the space discretisation is uniform.
\paragraph{Algorithm.}\label{algo} In the following, we state the algorithm  for computing the discrete values $\overline U_i^{n+1}$ at $t^{n+1}$. Given the initial data and boundary conditions and the number $CFL \leq 1$, we start with the known discrete averaged values $(\overline{ U}_i^n)$ for $i=0,...,N+1$ at $t^n$. As long as ($t<T$) one has to do:\\
\\
1) Computation of $\widetilde{A}_i$ for $i=0,...,N$ where $\widetilde{A}_i=A\left(\dfrac{\overline U_i^n+\overline U_{i+1}^n}{2}\right)$.\\
\\
2) Computation of $r_i^1$,$r_i^2$ and $\lambda_i^1$, $\lambda_i^2$ set respectively as the eigenvectors and eigenvalues of $\widetilde{A}_i$.\\
\\
3) Computation of $\Delta t$, such that $\dfrac{\Delta t}{\Delta x} \leq \dfrac{CFL}{a_{i+1/2}}$.\\
\\
4) Computation of ${\overline U_{i-1/2}^*}$ for $i=1,...,N+1$ by solving the linearized Riemann problem.\\
\\
In fact we have 3 cases:\\
\\
$\bullet$ if $\lambda_i^1$,$\lambda_i^2$<0 then ${\overline U_{i-1/2}^*}=\overline U_i^n$.\\
\\
$\bullet$ if $\lambda_i^1$,$\lambda_i^2$>0 then ${\overline U_{i-1/2}^*}=\overline U_{i-1}^n$.\\
\\
$\bullet$ if $\lambda_i^1<0$, $\lambda_i^2>0$ then for:\begin{equation}\left\{ \begin{array}{l}\dsp x<\lambda_i^1t \quad \mbox{one has} \quad {\overline U_{i-1/2}^*}=\overline U_{i-1}^n, \\ \\
        \dsp x>\lambda_i^1t \quad \mbox{or} \quad x<\lambda_i^2t  \quad \mbox{one has} \quad {\overline U_{i-1/2}^*}=\overline U_{i}^n-(R^{-1}[U])_2r_i^2=\overline U_{i-1}^n+(R^{-1}[U])_1r_i^1,\\ \\
        \dsp x>\lambda_i^2t \quad \mbox{one has} \quad {\overline U_{i-1/2}^*}=\overline U_{i}^n,\end{array} \right. \end{equation}
with $R=(r_i^1|r_i^2)$ and $[U]=\overline U^n_{i}-\overline U^n_{i-1}$.\\
\\
5) Computation of $F({\overline U_{i-1/2}^*})$.\\
\\
6) Computation of $\overline U^{n+1}_i=\overline U^n_i-\dfrac{\Delta t}{\Delta x}(F_{i+1/2}^n-F_{i-1/2}^n)$ for $i=1,...,N$.\\
\\
We repeat this algorithm for the new level of time $(t^{n+1}+\Delta t)$, until we reach the required final time $T$.
\subsubsection{Second order finite-volume scheme: MUSCL-RK2}\label{secMUSCLRK2}
A drawback of the Roe scheme (such as Godunov) is to be very diffusive. A remedy for this situation is through the extension of the scheme to the second order in space, associated to a second order Runge-Kutta scheme in time. In fact, we would like to reduce both numerical dissipation and dispersion within the hyperbolic component $S_1(.)$.
To this end, high order reconstructed
states at each interface have to be considered, following the classical MUSCL approach~\cite{Leer79} (Monotonic Upstream Scheme for Conservation Laws). To prevent the spurious oscillations that would occur around discontinuities or shocks, we suggest to use the ``minmod" limiter, designed to generate  slope limited, reconstructed left and right states for each cell that are used to calculate the flux at the interfaces. The implementation of this scheme is very easy and provides a natural extension of the Roe scheme described above. In fact, the main interests seem to be, after the tests, a gain of precision and stability. 

The steps of the second-order reconstruction are as follows: \\
1. Using the discrete averaged values $\overline U_i^n$, we construct the slopes $S_i^n$ using the ``minmod" limiter as the reconstruction must be non oscillatory in some sense, see~\cite{GR91}. We consider: 
\begin{equation}\label{minmod} 
S^{n}_i=\mbox{\emph{minmod}} \left(\dfrac{\overline U_{i+1}-\overline U_i}{x_{i+1}-x_i}, \dfrac{\overline U_{i}-\overline U_{i-1}}{x_{i}-x_{i-1}}\right) 
\end{equation}  where the function \emph{minmod} is defined on $\RR^2$ by
\begin{equation*}\mbox{\emph{minmod}}(a,b)=\left\{ \begin{array}{l}
        min (|a|,|b|) sgn(a) \quad \mbox{if} \quad sgn(a)=sgn(b),\\ \\
        0 \quad \mbox{else}.
    \end{array} \right. \end{equation*}
2. On the cell  $m_i = ]x_{i-1/2}, x_{i+1/2}[$, the solution is approached by:
\begin{equation*} U^{n}(x)= \overline U_i^{n} +(x-x_i)S_i^n. \end{equation*} 
3. We compute the numerical flux at the interfaces:
\begin{equation*}F^n_{i+1/2}=\widetilde{F}(\overline U_i^{n,+},\overline U_{i+1}^{n,-}) \quad \mbox{and} \quad F^n_{i-1/2}=\widetilde{F}(\overline U_{i-1}^{n,+},\overline U_{i}^{n,-})\end{equation*} with: 
\begin{equation*}\left\{ \begin{array}{l}
        \dsp \overline U_i^{n,+}\ =\ \overline U_i^{n}+\dfrac{\Delta x}{2}S_i^n\\ \\
        \overline U_{i+1}^{n,-}\ =\  \overline U_{i+1}^{n}-\dfrac{\Delta x}{2}S_{i+1}^n.
    \end{array} \right. \end{equation*}
4. We compute $\widetilde{U}_i^{n+1}$ by the application of~\eqref{volfini}, thus the scheme is given as follows:
\begin{equation*}\widetilde{U}^{n+1}_{i}= \overline U_{i}^n - \dfrac{\Delta t}{\Delta x}(F_{i+1/2}^n-F_{i-1/2}^n). \end{equation*}

As far as time discretization is concerned, we use the second-order explicit Runge–Kutta ``RK2" 
method which is described in the following. Given the ODE $\dfrac{dy}{dt}=f(t,y)$, one has,
\begin{eqnarray}\label{RK2}
    y^{n+1}&=&y^n+\dfrac{h}{2}\Big(f(t^n,y^n)+f(t^{n+1},\tilde{y}^{n+1})\Big),
\end{eqnarray}
with $\tilde{y}^{n+1}=y^n+hf(t^n,y^n)$, and $t^{n+1}=t^{n}+h$.\\
Applying~\eqref{RK2} to~\eqref{volfini}, we obtain the following modified scheme ``MUSCL-RK2":
\begin{equation}\label{RK2MUSCL}\overline U^{n+1}_i=\overline U^n_i-\dfrac{\Delta t}{2\Delta x}(F_{i+1/2}^n-F_{i-1/2}^n+F_{i+1/2}^{n+1}-F_{i-1/2}^{n+1}),\end{equation}
with \begin{equation*}F^{n+1}_{i+1/2}=\widetilde{F}(\widetilde{U}_i^{n+1,+},\widetilde{U}_{i+1}^{n+1,-}) \quad \mbox{and} \quad F^{n+1}_{i-1/2}=\widetilde{F}(\widetilde{U}_{i-1}^{n+1,+},\widetilde{U}_{i}^{n+1,-})  \ ,
\end{equation*}
$-\widetilde{F}$: the numerical flux determined as in the first order  VFRoe method, given by \eqref{numflux}.\\
$- \widetilde{U}^{n+1}_{i}$ is computed as follows: 
\begin{equation*}\widetilde{U}^{n+1}_{i}= \overline U_{i}^n - \dfrac{\Delta t}{\Delta x}(F_{i+1/2}^n-F_{i-1/2}^n) .
\end{equation*}
$-\widetilde{U}^{n+1,+}_{i}=\widetilde{U}^{n+1}_{i}+\dfrac{\Delta x}{2} \widetilde{S}_{i}^{n+1}.$ \\
$-\widetilde{U}^{n+1,-}_{i+1}=\widetilde{U}^{n+1}_{i+1}-\dfrac{\Delta x}{2} \widetilde{S}_{i+1}^{n+1}.$\\
$-\widetilde{S}_{i}^{n+1}$ and $\widetilde{S}_{i+1}^{n+1}$ are associated respectively to $\widetilde{U}_{i}^{n+1}$ and $\widetilde{U}_{i+1}^{n+1}$ by \eqref{minmod}.

\subsubsection{Higher order finite-volume scheme: WENO5-RK4}\label{WENO5sec}
The second order schemes are known to degenerate to first
order accuracy at smooth extrema. To reach higher order accuracy in smooth regions and a good resolution around discontinuities, we implement fifth-order accuracy WENO reconstruction, following~\cite{JiangShu96,Shu98}, where Jiang and Shu constructed third and fifth order finite difference WENO schemes in multi-space dimensions with a general framework for the design of the smoothness indicators and nonlinear weights. To automatically achieve high order accuracy and non-oscillatory property near discontinuities, WENO schemes use the idea of adaptive stencils in the reconstruction procedure based on the local smoothness of the numerical solution. We would like to mention also the previous studies~\cite{CLM,BCLMT,LannesMarche14}, where it is shown that for the study of dispersive waves, it is necessary to use high-order schemes to prevent the corruption of the dispersive properties of the model by some dispersive truncation errors linked to second-order schemes.
Using the same reconstruction proposed in~\cite{BCLMT}, we consider a cell $m_i$, and the corresponding constant averaged value $\overline{U}_i^n=(\overline{\zeta}_i^n,\overline{v}_i^n)$, with a constant space step $\Delta x$. This approach, provides high order reconstructed left and right values  $\overline{U}_i^{n,-}$ and $\overline{U}_i^{n,+}$, built following the five points stencil,
and introduced as follows:
\begin{equation}\label{recval}
    \overline{U}_i^{n,+}=\overline{U}_i^{n}+\dfrac{1}{2}\overline{\delta U}_i^{n,+} \quad \mbox{and} \quad \overline{U}_i^{n,-}=\overline{U}_i^{n}-\dfrac{1}{2}\overline{\delta U}_i^{n,-},
\end{equation}
where $\overline{\delta U}_i^{n,+}$ and $\overline{\delta U}_i^{n,-}$ are defined as follows:
\begin{eqnarray}\label{deltaUinplus}
    \overline{\delta U}_i^{n,+}&=\dfrac{2}{3}(\overline{U}_{i+1}^{n}-\overline{U}_{i}^{n})+\dfrac{1}{3}(\overline{U}_{i}^{n}-\overline{U}_{i-1}^{n})-\dfrac{1}{10}(-\overline{U}_{i-1}^{n}+3\overline{U}_{i}^{n}-3\overline{U}_{i+1}^{n}+\overline{U}_{i+2}^{n})\nn\\&\quad-\dfrac{1}{15}(-\overline{U}_{i-2}^{n}+3\overline{U}_{i-1}^{n}-3\overline{U}_{i}^{n}+\overline{U}_{i+1}^{n})
\end{eqnarray} 
\begin{eqnarray}\label{deltaUinmoins}
    \overline{\delta U}_i^{n,-}&=\dfrac{2}{3}(\overline{U}_{i}^{n}-\overline{U}_{i-1}^{n})+\dfrac{1}{3}(\overline{U}_{i+1}^{n}-\overline{U}_{i}^{n})-\dfrac{1}{10}(-\overline{U}_{i-2}^{n}+3\overline{U}_{i-1}^{n}-3\overline{U}_{i}^{n}+\overline{U}_{i+1}^{n})\nn\\&\quad-\dfrac{1}{15}(-\overline{U}_{i-1}^{n}+3\overline{U}_{i}^{n}-3\overline{U}_{i+1}^{n}+\overline{U}_{i+2}^{n})
\end{eqnarray}
and the coefficients $\dfrac{2}{3}$, $\dfrac{1}{3}$, $\dfrac{-1}{10}$ and $\dfrac{-1}{15}$ are set in order to obtain better dissipation and dispersion properties
in the truncature error. We consider the following modified scheme:
\begin{equation}\label{WENO5}\overline{U}^{n+1}_{i}= \overline U_{i}^n - \dfrac{\Delta t}{\Delta x}\Big(\widetilde{F}(\overline{U}_i^{n,+},\overline{U}_{i+1}^{n,-})-\widetilde{F}(\overline{U}_{i-1}^{n,+},\overline{U}_{i}^{n,-})\Big). \end{equation}
To reduce spurious oscillations near discontinuities, we apply the same limitation procedure as in~\cite{BCLMT}, preserving the scheme positivity and the high order accuracy. Thus scheme~\eqref{WENO5} becomes
\begin{equation}\label{WENO5LIM}\overline{U}^{n+1}_{i}= \overline U_{i}^n - \dfrac{\Delta t}{\Delta x}\Big(\widetilde{F}(^L\overline{U}_i^{n,+},^L\overline{U}_{i+1}^{n,-})-\widetilde{F}(^L\overline{U}_{i-1}^{n,+},^L\overline{U}_{i}^{n,-})\Big). \end{equation}
We define the limited high-order values as follows:
\begin{equation}\label{LIMrecval}
    ^L\overline{U}_i^{n,+}=\overline{U}_i^{n}+\dfrac{1}{2}L_{i}^{+}(\overline U^n) \quad \mbox{and} \quad ^L\overline{U}_i^{n,-}=\overline{U}_i^{n}-\dfrac{1}{2}L_{i}^{-}(\overline U^n).
\end{equation}
Using the following limiter,
\begin{equation*}L(u,v,w)=\left\{ \begin{array}{l}
        min (|u|,|v|,2|w|) \ sgn(u) \quad \mbox{if} \quad sgn(u)=sgn(v),\\ \\
        0 \quad \mbox{else},
    \end{array} \right. \end{equation*}
we define the limiting process as,
\begin{equation*}L_{i}^{+}(\overline U^n)=L(\delta \overline{U}_{i}^{n},\delta \overline{U}_{i+1}^{n},\overline{\delta U}_i^{n,+})\quad \mbox{and} \quad L_{i}^{-}(\overline U^n)=L(\delta \overline{U}_{i+1}^{n},\delta \overline{U}_{i}^{n},\overline{\delta U}_i^{n,-}),\end{equation*}
with $\delta \overline{U}_{i+1}^{n}=\overline U_{i+1}^{n}-\overline U_{i}^{n}$ and $\delta \overline{U}_{i}^{n}=\overline U_{i}^{n}-\overline U_{i-1}^{n}$ are upstream and downstream variations, and $\overline{\delta U}_i^{n,+} $ and $\overline{\delta U}_i^{n,-}$ taken from~\eqref{deltaUinplus} and~\eqref{deltaUinmoins}.

The limited high order reconstructions stated above must be performed on both conservative variables $\overline{U}_i^n=(\overline{\zeta}_i^n,\overline{v}_i^n)$. We would like to mention that the resulting finite volume scheme preserve motionless steady states, $\zeta=\text{cst}$ and $v=0$.

As far as time discretization is concerned, we use the fourth-order explicit Runge–Kutta ``RK4" method which is described in the following. Given the ODE $\dfrac{dy}{dt}=f(t,y)$, one has,
\begin{eqnarray}\label{RK4}
    k_1&=&f(t^n,y^n)\nn,\\
    k_2&=&f(t^n+\dfrac{h}{2},y^n+h\dfrac{k_1}{2})\nn,\\
    k_3&=&f(t^n+\dfrac{h}{2},y^n+h\dfrac{k_2}{2})\nn,\\
    k_4&=&f(t^n+h,y^n+hk_3)\nn,\\
    y^{n+1}&=&y^n+\dfrac{h}{6}\big(k_1+2k_2+2k_3+k_4\big),
\end{eqnarray}
with $t^{n+1}=t^n+h$. Applying~\eqref{RK4} to~\eqref{WENO5LIM}, one gets the ``WENO5-RK4" scheme.
\subsection{Finite difference scheme for the dispersive part}\label{Dispersive}
First of all, let us detail how to construct nodal values of the unknowns (which are the ones used for a finite difference discretization) from the
cell-averaged value computed by a finite volume scheme and vice versa.

We denote by $U_i^n$ the nodal value of $U$ at the $i$th node $(x_{i+1/2})_{i\in[0,N]}$ and at time $t_n$ i.e. $U_i^n$  is an approximation of $U(x_{i+1/2},t^n)$.
The finite volume-finite difference mix  imply to switch between the cell-averaged and nodal values for each unknown and at each time step.
To this end, we use the fifth-order accuracy WENO reconstruction, that allows to approximate the nodal values (i.e finite difference unknowns) 
$(U_i^n)_{i =1,N+1}$ in terms of the cell-averaged values (i.e finite volume unknowns) $(\overline U_i^n)_{i =1,N}$.   The relation is given by:
    \begin{equation}\label{switchFVDF}
    U_i^n=\dfrac{1}{30}\overline U_{i-2}^n-\dfrac{13}{60}\overline U_{i-1}^n+\dfrac{47}{60}\overline U_i^n+\dfrac{9}{20}\overline U_{i+1}^n-\dfrac{1}{20}\overline U_{i+2}^n +\OO(\Delta x^5), \quad 1\leq i \leq N+1,
    \end{equation}
 with adaptations at the boundaries following the method presented in Section~\ref{BCsec}. 
 One can easily recover the relation that allows to determine the cell-averaged values $(\overline U_i^n)_{i \in [1:N]}$ in terms of the nodal values $(U_i^n)_{i \in [1:N+1]}$.
    
We can easily check that~\eqref{switchFVDF} preserve the steady state at rest and that this formula is precise up to order $\OO(\Delta x^5)$
terms, thus preserving the global order of the scheme.

We can now proceed by explaining how we compute the solution operator $S_2(.)$ associated to the dispersive part of the equations. 
Let us recall the system~\eqref{disp} corresponding to the operator $S_2(.)$, given in Section~\eqref{SSsec}.
\begin{equation}\label{disp1}\left\{ \begin{array}{l}
        \dsp \partial_{ t}\zeta \ =\  0,\\ \\
        \dsp  \partial_{ t}  v-\dfrac{1}{\alpha}\big((\gamma+\delta)\partial_x \zeta +\epsilon\partial_x(q_3(\epsilon\zeta)  {v}^2)\big)\\+(I+\mu\nu\alpha T[0])^{-1}\Big[\dfrac{1}{\alpha}\big((\gamma+\delta)\partial_x \zeta +\epsilon\partial_x(q_3(\epsilon\zeta)  {v}^2)\big)+\mu\epsilon Q_1(v)+\mu\nu Q_2(\zeta)+\mu\epsilon\nu Q_3(\zeta)\Big]=0.
    \end{array} \right. \end{equation}
where the operators $Q_1$, $Q_2$ and $Q_3$ are explicitly given in~\eqref{Q1},~\eqref{Q2} and~\eqref{Q3}.
For the sake of simplicity, we detail the numerical resolution of \eqref{disp1} using an explicit Euler in time scheme. Standard extensions to second 
and fourth order Runge-Kutta method has been used according to the order of the space derivative operators as done in the previous section.

The finite discretization of the system~\eqref{disp1} leads to the following discrete problem:
\begin{equation}\label{disp1disc}\left\{ \begin{array}{l}
        \dsp \dfrac{\zeta^{n+1}-\zeta^n}{\Delta t} \ =\  0,\\ \\
        \dsp \dfrac{v^{n+1}-v^n}{\Delta t}-\dfrac{1}{\alpha}(\gamma+\delta)D_1 (\zeta^{n})
        -2\dfrac{\epsilon}{\alpha}q_3(\epsilon\zeta^{n})v^{n}D_1(v^{n})-\dfrac{\epsilon^2}{\alpha}q'_3(\epsilon\zeta^{n})D_1(\zeta^{n})(v^{n})(v^{n})
        \\+(I-\mu\nu\alpha D_2)^{-1} \Big[ \dfrac{1}{\alpha}(\gamma+\delta)D_1(\zeta^{n})+2\dfrac{\epsilon}{\alpha}q_3(\epsilon\zeta^{n})v^{n}D_1(v^{n})+\dfrac{\epsilon^2}{\alpha}q'_3(\epsilon\zeta^{n})D_1(\zeta^{n})(v^{n})(v^{n})\\ \medskip+\mu\epsilon Q_1(v^n)+\mu\nu\epsilon Q_2(\zeta^n)+\mu\epsilon\nu Q_3(\zeta^n)\Big]=0,
    \end{array} \right. \end{equation}
with
\begin{equation*}\label{discQ1}Q_1(v^n)=2 \kappa D_1(v^{n})D_2(v^{n}),\end{equation*}
\begin{equation*}\label{discQ2}Q_2(\zeta^n)=\kappa_2D_1\Big[\zeta^{n}D_1\Big((I-\mu\nu\alpha D_2)^{-1}(\gamma+\delta)D_1(\zeta^{n})
    \Big)\Big],
\end{equation*}
\begin{equation*}\label{discQ3}Q_3(\zeta^n)=-\kappa_1\zeta^{n}D_2\Big[(I-\mu\nu \alpha D_2)^{-1}(\gamma+\delta)D_1(\zeta^{n})
    \Big].\end{equation*}
The system~\eqref{disp1disc} is solved at each time step using a classical finite-difference technique, where the matrices $D_1$ and $D_2$ are the classical centered discretizations of the derivatives $\partial_x$ and $\partial^2_x$ given below. 

The first formula is the second-order formula called ``DF2", where the spatial derivatives are given as follows:
\begin{equation*}\label{Diff1}
    (\partial_x U)_i=\dfrac{1}{2\Delta x}(U_{i+1}-U_{i-1}),
\end{equation*}
\begin{equation*}\label{Diff2}
    (\partial_x^2 U)_i=\dfrac{1}{\Delta x^2}(U_{i+1}-2U_{i}+U_{i-1}).
\end{equation*}
For time discretization, the second-order formula ``DF2" is associated to a second-order classical Runge-Kutta ``RK2" scheme, and thus one obtains the ``DF2-RK2" scheme.

The second formula is the fourth-order formula called ``DF4" where the spatial derivatives are given as follows:
\begin{equation*}\label{Diff1-4}
    (\partial_x U)_i=\dfrac{1}{12\Delta x}(-U_{i+2}+8U_{i+1}-8U_{i-1}+U_{i-2}),
\end{equation*}
\begin{equation*}\label{Diff2-4}
    (\partial_x^2 U)_i=\dfrac{1}{12\Delta x^2}(-U_{i+2}+16U_{i+1}-30U_{i}+16U_{i-1}-U_{i-2}).
\end{equation*}
For time discretization, the fourth-order formula ``DF4" is associated to a fourth-order classical Runge-Kutta ``RK4" scheme, and thus one obtains the ``DF4-RK4" scheme.
\subsection{ Boundary conditions}\label{BCsec}
In the following section, we show how to treat the boundary conditions for the hyperbolic and dispersive parts of the splitting scheme. 
Suitable relations are imposed on both cell-averaged and nodal quantities. We only treat either periodic boundary conditions or 
reflective boundary conditions. We detail now how we have implemented these boundary conditions for the hyperbolic part and the dispersive part
of the numerical scheme.

For the hyperbolic part, we have introduced ‘‘ghosts cells’’  respectively at the upstream and downstream boundaries of the domain. 
The imposed relations on the cell-averaged quantities  are the following:

$\bullet$ $\overline U_{-k+1}=\overline U_{N-k+1}$, and $\overline U_{N+k}=\overline U_{k}$, $k \geq 1$, for periodic conditions on upstream and downstream boundaries.

$\bullet$ $\overline \zeta_{-k+1}=\overline \zeta_{-k}$, $\overline v_{-k+1}=-\overline v_{-k}$ and $\overline \zeta_{N+k-1}=\overline \zeta_{N+k}$, $\overline v_{N+k-1}=-\overline v_{N+k}$, $k \geq 1$, for reflective conditions on the left and right boundaries.\\

For the dispersive part of the splitting, the boundary conditions are simply imposed on the nodal values that are located outside of the domain, in order to maintain centered formula at the boundaries, while keeping a regular structure in the discretized model:

$\bullet$ $U_{-k+1}= U_{N-k+1}$, and $U_{N+k}=U_{k}$, $k \geq 1$, for periodic conditions on upstream and downstream boundaries.

$\bullet$ $\zeta_{-k+1}= \zeta_{-k}$, $ v_{-k+1}=-v_{-k}$ and $ \zeta_{N+k-1}=\zeta_{N+k}$, $ v_{N+k-1}=-v_{N+k}$, $k \geq 1$, for reflective conditions on upstream and downstream boundaries.

\section{Numerical validations}\label{NVsec}
In this section, several numerical tests are performed in both one and two layers cases in order to validate the numerical efficiency and accuracy of the improved Green-Naghdi model~\eqref{GNCH6dim}. We first consider several numerical tests in the one layer case without any surface tension i.e $\bo^{-1}=0$. We begin by studying the propagation of a solitary wave over a flat bottom. We compare our numerical solution with an analytic one (up to an $\OO(\mu^2)$ remainder) at several times and for different orders of discretizations and show that our numerical scheme is very efficient and accurate. To evaluate the influence of the nonlinear and dispersive terms we study the collision between two solitary waves traveling in opposite directions (head-on collision). We then study the breaking of a Gaussian hump into two solitary waves. Finally, we study the dam-break problem supplemented by a comparison between the second and fifth order accuracy in order to show the ability of the higher order numerical scheme in dealing with discontinuities. We used the value $\alpha=1$ in the aforementioned cases. In fact, we obtained very similar results when performing the same simulations with $\alpha=1.159$. Secondly, two numerical simulations are performed in the two layers cases. In the first one, we compare our results with numerical data from~\cite{DucheneIsrawiTalhouk16}, where we show that a very good matching  is observable if $\alpha$ is carefully chosen as in Section~\ref{Secalphachoice}. We then consider the dam-break problem in the two layers cases, where we test the ability of the splitting scheme to compute dispersive shock waves with high accuracy. We would like to mention that in all the numerical tests, we use the WENO5 reconstruction for the hyperbolic part of the splitting scheme and a fourth order finite difference scheme ``DF4" for the dispersive part, both associated to a fourth-order classical Runge-Kutta ``RK4" time scheme. In every numerical simulation presented in the following section, we choose to use a CFL number equal to 1 in the algorithm stated in page~\pageref{algo}, in order to obtain a stable numerical scheme.
\subsection{Numerical validations in the one layer case}\label{NV1sec}
\subsubsection{Propagation of a solitary wave}~\label{PSWsec}
\begin{figure}[H]
    \centering
{\includegraphics[scale=0.44]{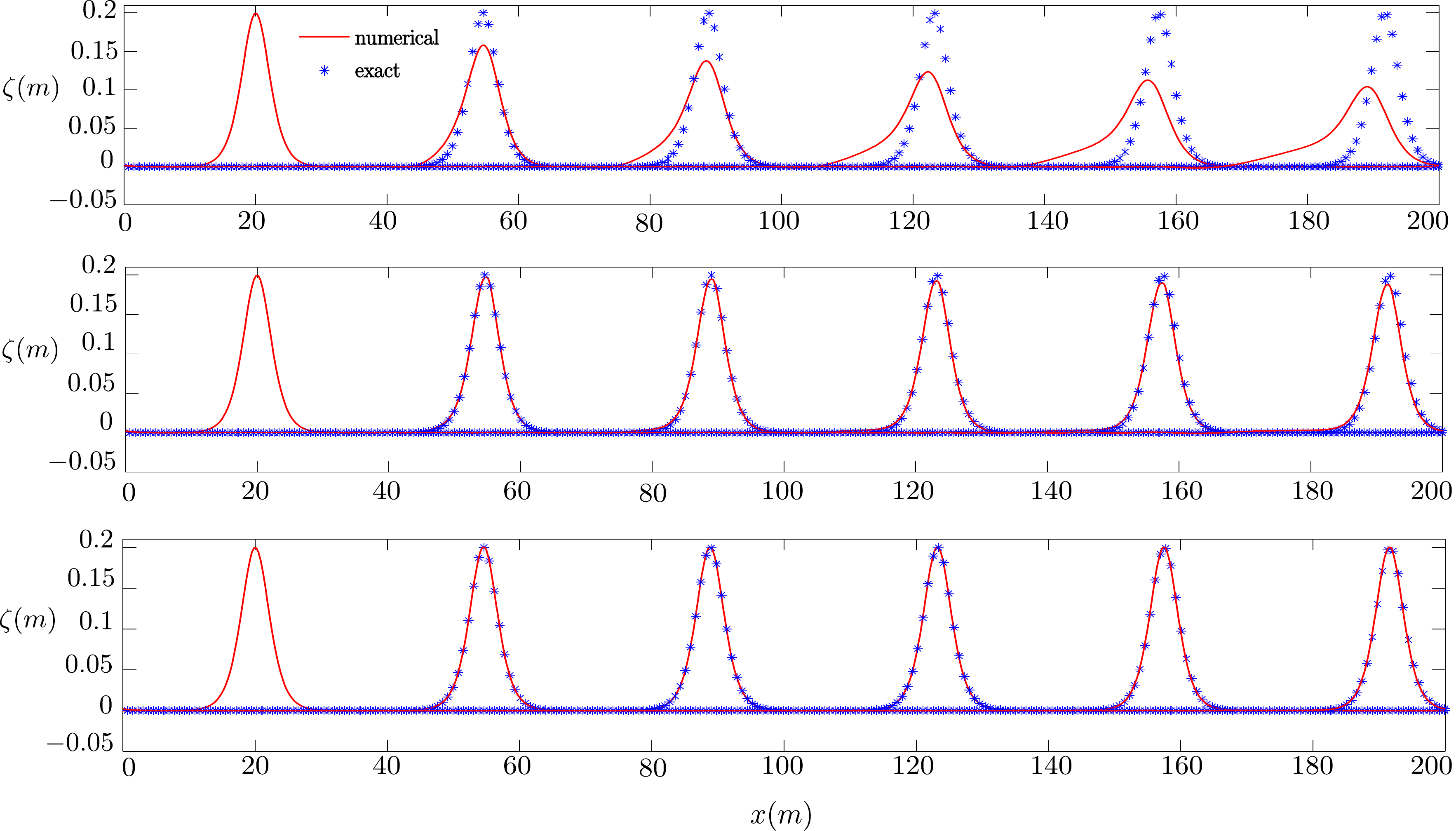}}
    \caption{Propagation of a solitary wave over a flat bottom: water surface profiles at t=0, 10, 20, 30, 40 and 50s. Top: FV1-DF2-Euler, middle: MUSCL-DF2-RK2 and bottom: WENO5-DF4-RK4. }
    \label{125}
\end{figure}

Here, we test the accuracy of our numerical scheme~\eqref{GNCH6dim} with $\alpha=1$, by using the exact solitary wave solutions of the one layer Green-Naghdi equations in the one-dimensional setting over a flat bottom (see~\cite{LannesMarche14}), given in variables with dimensions, by
\begin{equation}\label{soliton}\left\{ \begin{array}{l}
        \dsp \zeta(t,x)=a \ \text{sech}^2(k(x-ct))  ,\\ \\
        \dsp  v(t,x)=c \Big(\dfrac{\zeta(t,x)}{d_2+\zeta(t,x)}\Big),\\ \\
        \dsp  k=\dfrac{\sqrt{3 a}}{2d_2\sqrt{d_2+a}}, \quad c=\sqrt{g(d_2+a)},
    \end{array} \right. \end{equation}
    where we recall that $d_2$ is depth of the fluid when considering the one layer case.
Such solitary waves are also solution of the improved Green-Naghdi model~\eqref{GNCH6dim} up to an $\OO(\mu^2)$ remainder. 
We consider the propagation of a solitary wave initially centered at $x_0=20 \ m$, of relative amplitude $a=0.2 \ m$, over a constant water depth $d_2=1 \ m$. The computational domain is $200 \ m$ in length and discretized with 1280 cells. The single solitary wave propagates from left to right. In this test, since the solitary wave is initially far from boundaries, the boundary conditions do not affect the computation, thus we choose to impose periodic boundary conditions on each boundary for the sake of simplicity.
We compare the water surface profile of our numerical solution provided by the model~\eqref{GNCH6dim} (after setting the parameters $\gamma=0$ and $\delta=1$ corresponding to the one layer case), with the exact one given by~\eqref{soliton} at several times using the first, second and fifth order discretizations (see Figure~\ref{125}). One can easily remark that the fifth order discretization ``WENO5-DF4-RK4" provides very accurate solutions and reduces both numerical dissipation and dispersion, contrarily to the first order discretization ``FV1-DF2-Euler" which seems to be very diffusive. Indeed, looking at the amplitude and shape of the solitary wave at $t=50 \ s$ in the bottom of Figure~\ref{125}, we can observe an excellent agreement between numerical and exact solutions, unlike the second order discretization ``MUSCL-DF2-RK2" (middle of Figure~\ref{125}), generating a less accurate numerical solution. The preservation of the amplitude and shape of the solitary wave computed using the fifth order scheme, even after $200 \ m$ indicate that the governing equations have been accurately discretized in space and time.

In what follows, we quantify the numerical accuracy of our numerical scheme by computing the numerical solution for this particular test case for an increasing number of cells $N$, over a duration $T=5\ s$. Starting with $N=80$ number of cells, we successively multiply the number of cells by two. The relative errors  $E_{L^2}(\zeta)$ and $E_{L^2}(v)$ on the water surface deformation and the averaged velocity presented in Table~\ref{L2err} are computed at $t= 5 \ s$, using the discrete $L^2$ norm $\|.\|_{2}$:
\begin{equation}
E_{L^2}(\zeta)=\dfrac{\|\zeta_{num}-\zeta_{sol}\|_{2}}{\|\zeta_{sol}\|_{2}}; \qquad E_{L^2}(v)=\dfrac{\|v_{num}-v_{sol}\|_{2}}{\|v_{sol}\|_{2}},
\end{equation}
where $(\zeta_{num}, v_{num})$ are the numerical solutions and $(\zeta_{sol}, v_{sol})$ are the analytical ones coming from~\eqref{soliton}. Very accurate results are thus obtained as an evaluation of the capacity of our
numerical method to compute in a stable way the propagation of a solitary wave.
\begin{table}[H]
\centering
 \begin{tabular}{  p{4cm} p{4cm}  p{3cm} }
FV1-DF2-Euler &MUSCL-DF2-RK2  & WENO5-DF4-RK4 
 \end{tabular}
 \begin{tabular}{ | c | c | c | c | c | c | c | }
\hline
N & $E_{L^2}(\zeta)$&$E_{L^2}(v)$& $E_{L^2}(\zeta)$&$E_{L^2}(v)$& $E_{L^2}(\zeta)$&$E_{L^2}(v)$\\
  \hline		
 80 & $5.79\times 10^{-1}$ & $5.56\times 10^{-1}$ & $5.57\times 10^{-1}$ & $5.30\times 10^{-1}$& $4.32\times 10^{-1}$ & $4.02\times 10^{-1}$\\
  160 & $4.30\times 10^{-1}$ & $4.07\times 10^{-1}$ & $3.54\times 10^{-1}$ & $3.27\times 10^{-1}$& $1.94\times 10^{-1}$ & $1.67\times 10^{-1}$\\
  320 & $3.04\times 10^{-1}$ & $2.83\times 10^{-1}$ & $1.76\times 10^{-1}$ & $1.54\times 10^{-1}$& $6.45\times 10^{-2}$ & $5.25\times 10^{-2}$ \\
 640 & $1.95\times 10^{-1}$ & $1.79\times 10^{-1}$ & $5.96\times 10^{-2}$ & $5.00\times 10^{-2}$& $1.16\times 10^{-2}$ & $9.30\times 10^{-3}$  \\
 1280& $1.14\times 10^{-1}$ & $1.04\times 10^{-1}$ & $1.38\times 10^{-2}$& $1.20\times 10^{-2}$ & $3.60\times 10^{-3}$& $3.40\times 10^{-3}$\\
  \hline  
\end{tabular}
 \caption{Propagation of a solitary wave over a flat bottom: relative $L^2$-error table for the conservative variables.}
 \label{L2err}
\end{table}
 \begin{remark}\label{remfilter}
We believe that the main reason for not obtaining the predicted order in each space discretization might be due to the fact that the analytic solution given in~\eqref{soliton} satisfies the model~\eqref{GNCH6dim} up to an $\OO(\mu^2)$ remainder, that is to say it is an approximate solution. A remedy for this situation could be through an ``iterative cleaning" technique that acts to damp the high frequency oscillations (i.e the oscillatory dispersive tails) that appears due to the remainder term of size $\OO(\mu^2)$. This technique has been used by several authors, see for instance~\cite{BC98,BDM07,ND08}. In this paper, we do not try to give some optimal convergence result and the filtering technique is left to future work.
\end{remark}
%
%
%
\subsubsection{Head-on collision of counter-propagating waves}\label{HOCsec}
We now investigate the interaction of solitary waves which allows us to evaluate the impact of non linearities and dispersive terms. To this end, we study the head-on collision of two solitary waves traveling in opposite directions. The initial data for the two counter-propagating solitary waves are given in~\eqref{soliton}. Many authors have set different models and numerical methods in order to numerically study this problem (see~\cite{CGHHS06,Arnaud14,MID14}). Unlike solitary waves in integrable systems, the one for the \emph {full Euler} equations are often followed by a non zero residual wave after interactions. The resulting residual has the form of dispersive trailing waves of small amplitude. In this test, we study the interaction of two solitary waves of equal amplitudes propagating in an opposite direction, initially located at $x=100 \ m$ and $x= 300 \ m$. The spatial domain is $400 \ m$ in length discretized using $1200$ cells. Periodic conditions are imposed on each boundary.
\begin{figure}[H]
    \centering
{\includegraphics[scale=0.23]{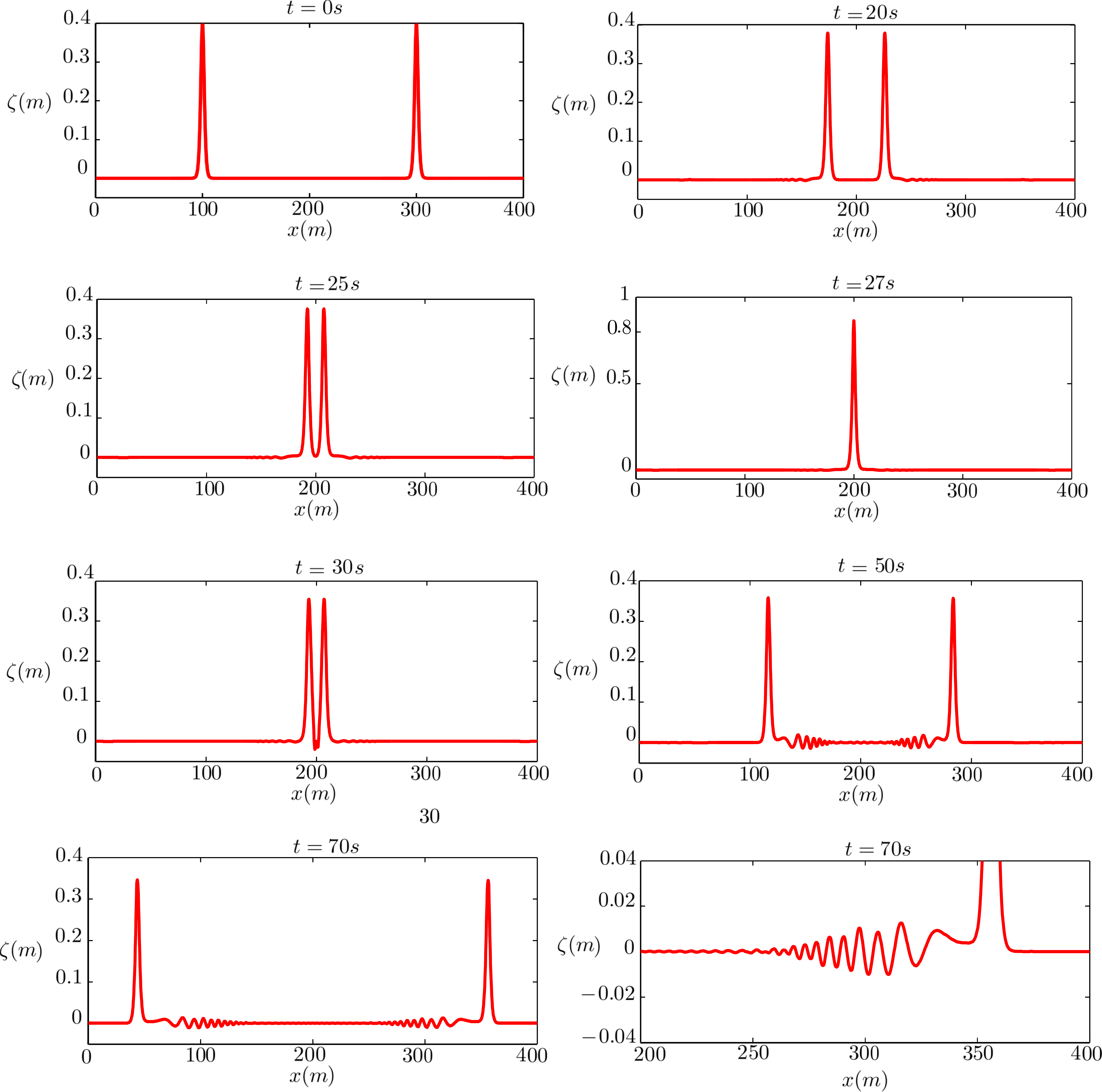}}
    \caption{Head on collision of two counter-propagating solitary waves: water surface profiles at several times during the propagation.}
    \label{HOC}
\end{figure}
The results are shown in Figure~\ref{HOC} at different propagation times. Before the collision, at times $t=20 \ s$ and $t= 25 \ s$ one can observe two dispersive tails of very small amplitude located to the left and right of each solitary wave. The generation of such dispersive tails is due to the $\OO(\mu^2)$ remainder term as mentioned in Remark~\ref{remfilter}. As expected, the waves collide to reach a maximum height larger than the sum of the
amplitudes of the two incident solitary waves at the time $t = 27 \ s$. After the collision, dispersive tails with small amplitudes appear clearly when zooming at $t= 70 \ s$, illustrating an appropriate characterization of the nonlinear interactions. Capturing this dynamics validate the high precision of our numerical scheme. The head-on collision is also studied in~\cite{MID14,Arnaud14}, leading to similar observations. 
\subsubsection{Breaking of a Gaussian hump}\label{HOWsec}
\begin{figure}[H]
    \centering
{\includegraphics[scale=0.23]{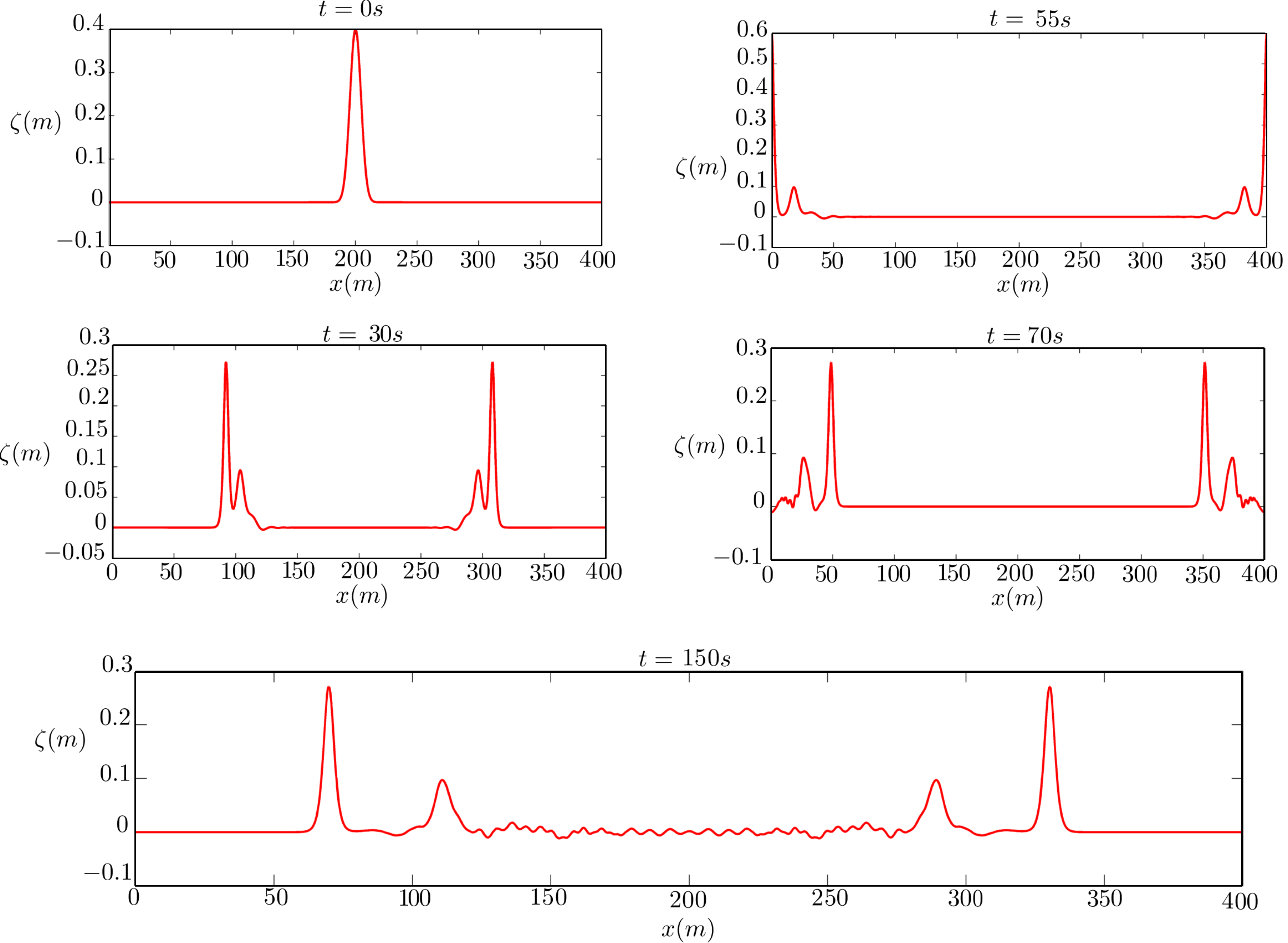}}
    \caption{Breakup of a Gaussian hump into two solitary waves traveling in opposite directions and dispersive tails. }
    \label{HOW}
\end{figure}

In this test, we consider the following initial data representing a heap of water,
\begin{equation*}
\zeta(0,x)=a e^{-\frac{1}{\lambda}(x-\frac{L}{2})^2}, \quad v(0,x)=0,
\end{equation*} 
where $a$ represents the amplitude, $\lambda$ represents the width and $L$ represents the length of the domain. Figure~\ref{HOW}, shows the overall behavior of the solutions using $a=0.4 \ m$ and $\lambda=40 \ m$ and $L=400 \ m$ discretized with $2000$ cells using periodic boundary conditions. The initial hump breaks up into two large solitary waves and smaller dispersive tails. These waves and the dispersive tails travel in opposite directions. 
 The results shown in Figure~\ref{HOW}, tends to confirm the ability of our numerical scheme to capture this dynamics accurately. Indeed, this typical behaviour is also studied in~\cite{MID14}, leading to very similar observations.
\subsubsection{Dam-break problem in the one layer case}\label{DB1sec}
We consider now a dam-break problem in the one layer case in order to test the ability of our numerical scheme to deal with discontinuities.
In general, discontinuous initial data of this type generates dispersive shock waves due to the dispersive effects~\cite{MGH10}. Analytic and computational studies of the dispersive shock waves in fully-nonlinear dispersive shallow water systems were carried out in~\cite{EGS06,MGH10}.
We would like to mention also the previous studies~\cite{CLM,BCLMT,LannesMarche14}, where it is shown that for the study of dispersive waves, it is necessary to use high-order schemes to prevent the corruption of the dispersive properties of the model by some dispersive truncation errors linked to second-order schemes. Indeed, this test is supplemented by a comparison between the second and fifth order accuracy exhibiting the ability of higher order schemes to capture the rapid oscillations in dispersive shock waves. We use the following initial data:
\begin{equation}\label{daminitialdata}
\zeta(0,x)=a[1+\tanh(250-|x|)], \quad v(0,x)=0,
\end{equation} 
where $a=0.2091 \ m$. The computational domain is the interval $x\in (-700,700)$ and discretized using 2800 cells. We choose to impose periodic boundary conditions.
\begin{figure}[H]\
    \centering
{\includegraphics[scale=0.52]{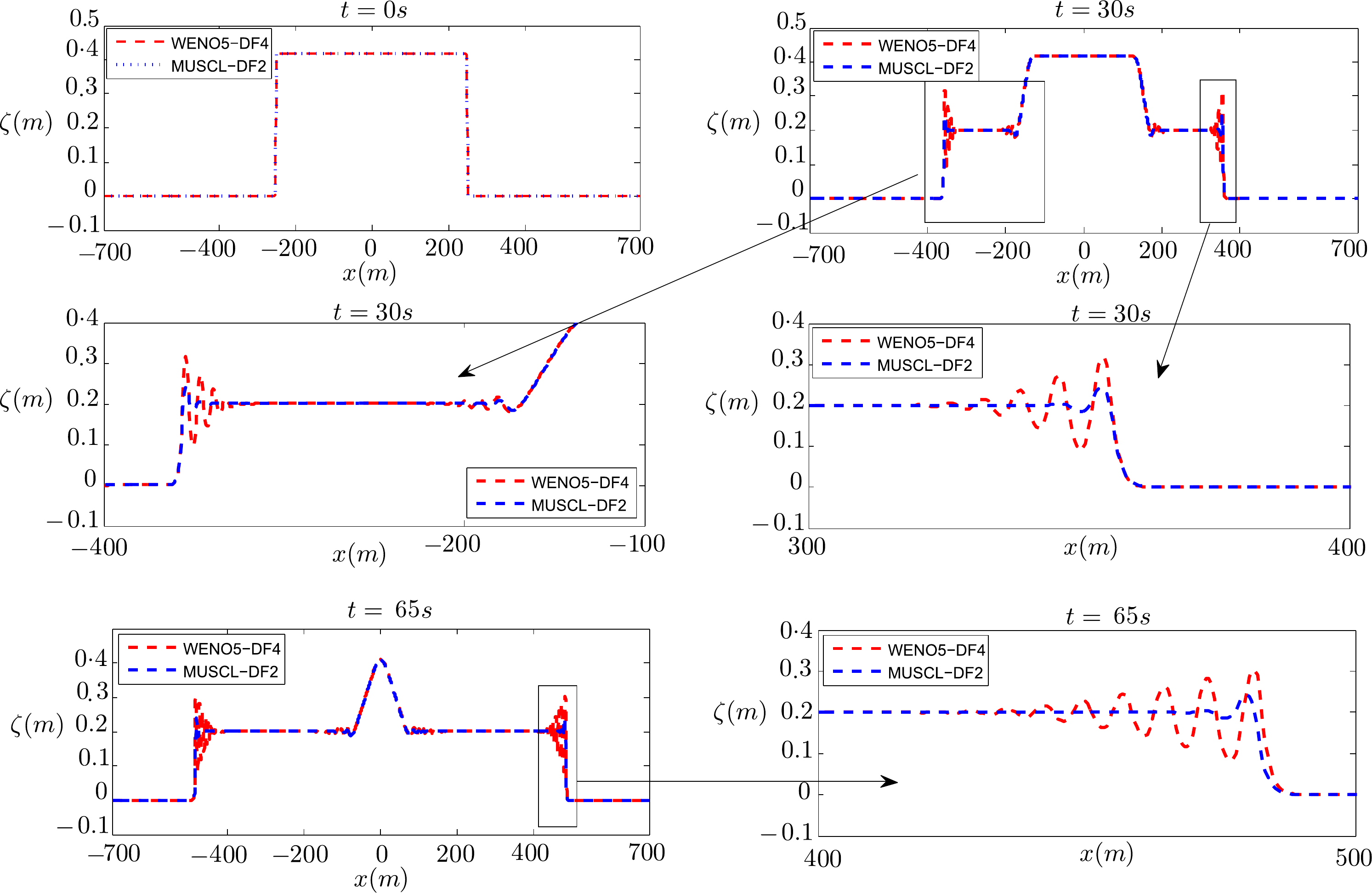}}
    \caption{Dam break in the one layer case: water surface profiles at several times}
    \label{dam1}
\end{figure}
    Figure~\ref{dam1} shows the results of the numerical simulations with two different orders of discretization, ``MUSCL-DF2-RK2" and ``WENO5-DF4-RK4", generating two dispersive tails propagating in opposite direction, on the left and right side of the “dam”, and two rarefaction waves that travel towards the center. A zoom on the shock and rarefaction waves at $t=30 \ s$ and $t=65 \ s$, shows clearly the corruption of the dispersive properties by the second order scheme. The dam-break problem was also studied in~\cite{EGS06,MGH10} and seems to fit well with our numerical observations. The numerical model proposed in~\cite{MID14} computes using a finite element method all the nonlinear dispersive terms, in particular the third order ones. These third order derivatives are present in our model but in order to improve the frequency dispersion we have proposed to factorize these terms, making it possible not to compute them. This formulation has the inconvenience of numerical diffusion. This is the reason why the dispersive tails observed in the dam-break problem in~\cite{MID14} have larger amplitude oscillations.
\subsection{Numerical validations in the two layers case}\label{NV2sec}
\subsubsection{Kelvin-Helmholtz instabilities}\label{KH}
In this case, we would like to highlight the importance of the choice of the parameter $\alpha$ in order to improve the frequency dispersion of the model~\eqref{GNCH6}, through the simulation of a sufficiently regular initial wave, following the numerical experiments performed in~\cite{DucheneIsrawiTalhouk16}.  In the aforementioned paper they introduce a new class of Green-Naghdi type models for the propagation of internal waves with improved frequency dispersion in order to prevent high-frequency Kelvin-Helmholtz instabilities. These models are obtained by regularizing the original Green-Naghdi one by slightly modifying the dispersion components using a class of Fourier multipliers. They represent three different choices of the Fourier multipliers, each one yields to a specific Green-Naghi model which they denote as follows:\\
\\
$\bullet$  ``original" as the classical Green-Naghi model introduced in~\cite{ChoiCamassa99}.\\
\\
 $\bullet$  ``regularized" which is a well-posed system for sufficiently small and regular data, even in absence of surface tension. In addition its dispersion relation fit the one of \emph {full Euler system} at order $\OO(\mu^3)$.\\
\\
$\bullet$   ``improved" whose dispersion relation is the same as the one of the \emph {full Euler system}.\\
\\
Using the spectral methods~\cite{Trefethen00} for the space discretization and the Matlab solver ode45, which is based on the fourth and fifth order Runge-Kutta-Merson method for time evolution, they numerically compute several of their Green-Naghdi systems. Several computations are made, with and without surface tension in order to observe how the different frequency dispersion may affect the appearance of Kelvin-Helmholtz instabilities. In order to compare with the numerical experiments done in~\cite{DucheneIsrawiTalhouk16} we choose the initial data $\zeta(0,x)=-e^{-4|x|^2}$ and $v(0,x)=0$ (represented by the dashed lines). The dimensionless parameters are set as follows: $\mu=0.1$, $\epsilon=0.5$, $\delta=0.5$, $\gamma=0.95$ and $\bo^{-1}=5\times10^{-5}$. The computational domain is the interval $x \in (-4,4)$ discretized with 512 cells using periodic boundary conditions. In all the following simulations we compute our numerical solution using the fifth order accuracy scheme ``WENO5-DF4-RK4".
  \begin{figure}[H]
    \centering
{\includegraphics[scale=0.68]{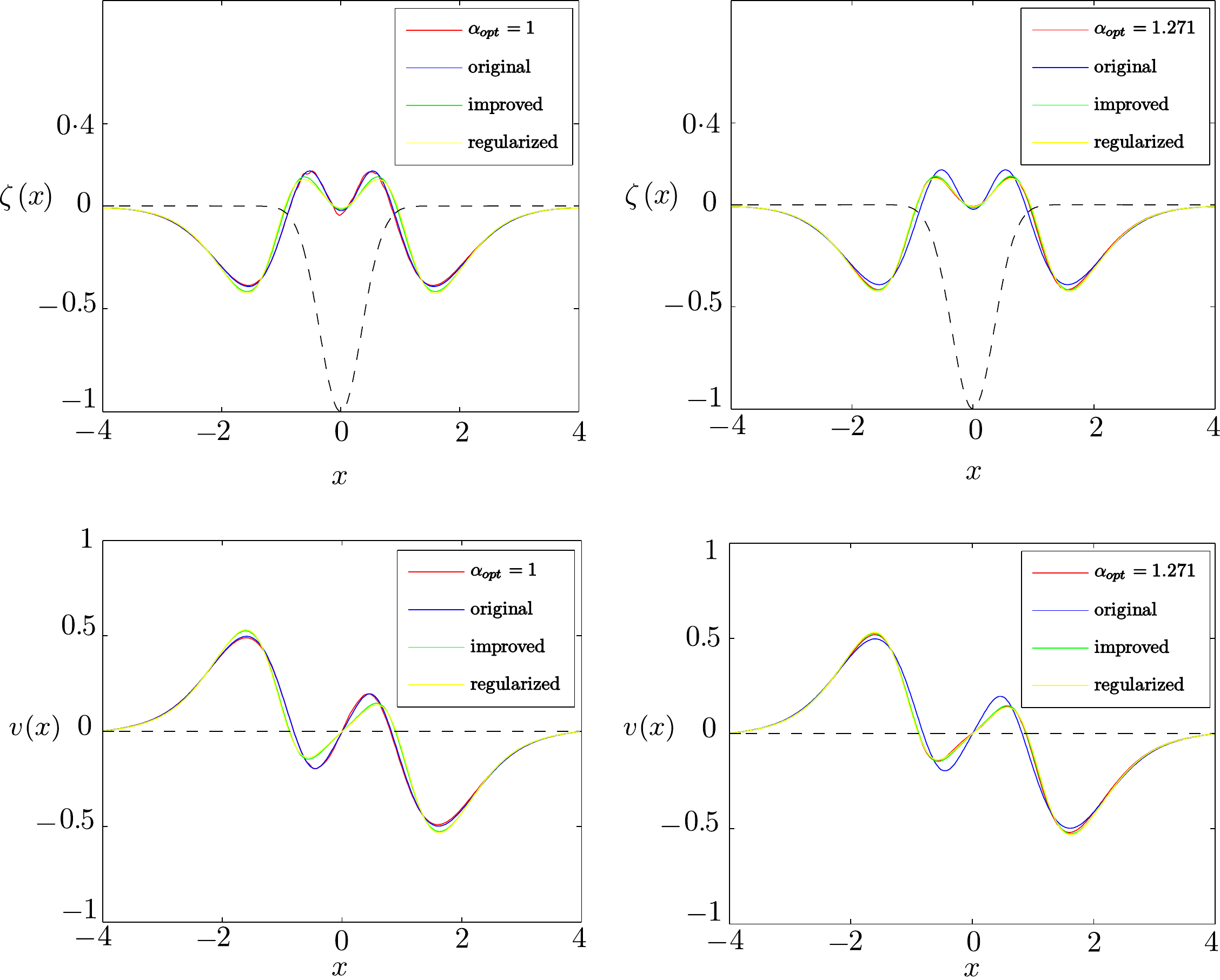}}
    \caption{Comparison with the Green-Naghdi models, with surface tension, at time $t=2$, for $\alpha_{opt}=1$ (left) and $\alpha_{opt}=1.271$ (right)}
    \label{DIT1}
\end{figure}
    \hspace{7cm}
\begin{figure}[H]
    \centering
{\includegraphics[scale=0.68]{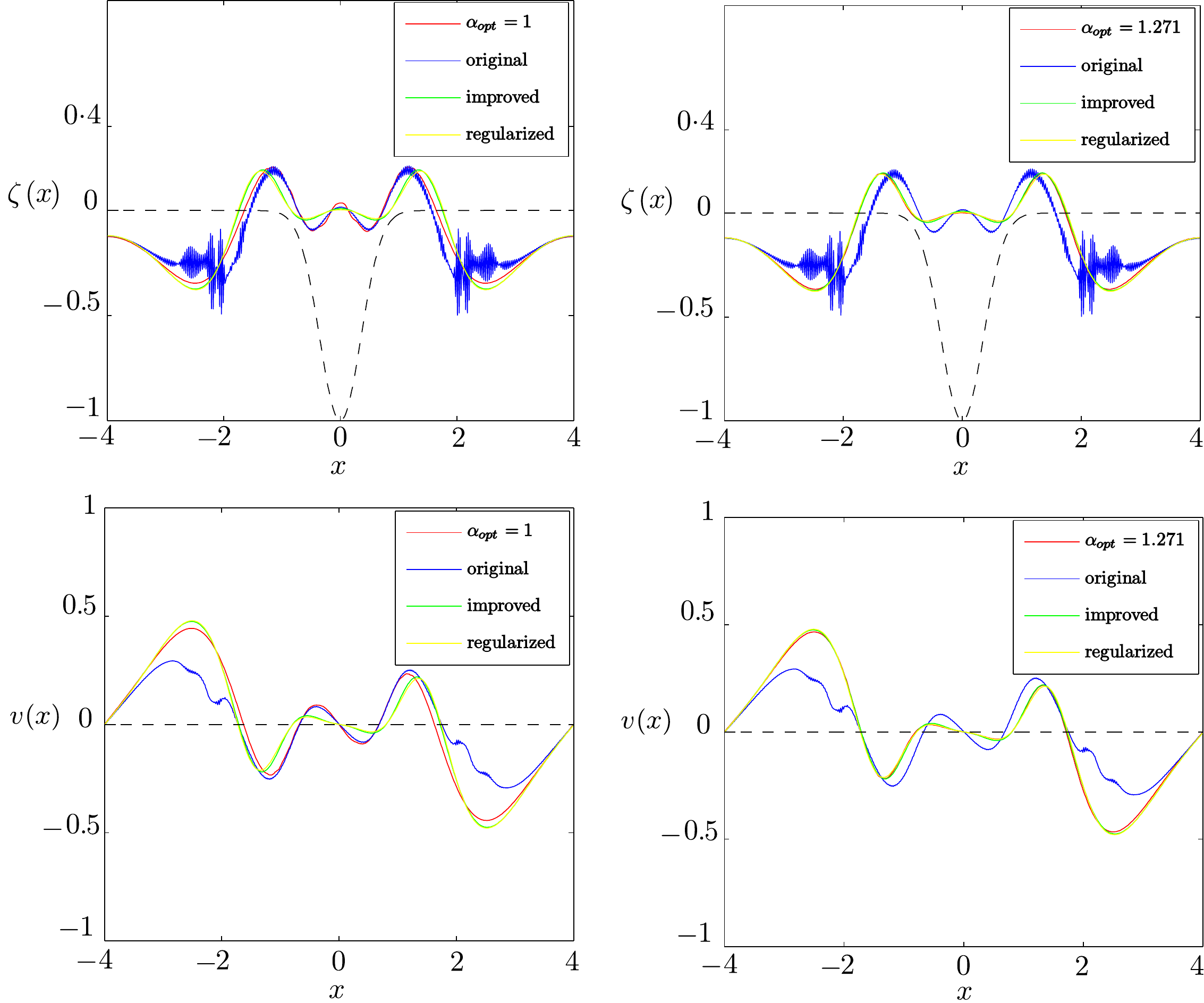}}
    \caption{Comparison with the Green-Naghdi models, with surface tension, at time $t=3$, for $\alpha_{opt}=1$ (left) and $\alpha_{opt}=1.271$ (right)}
    \label{DIT2}
\end{figure}
    \begin{figure}[H]
    \centering
{\includegraphics[scale=0.68]{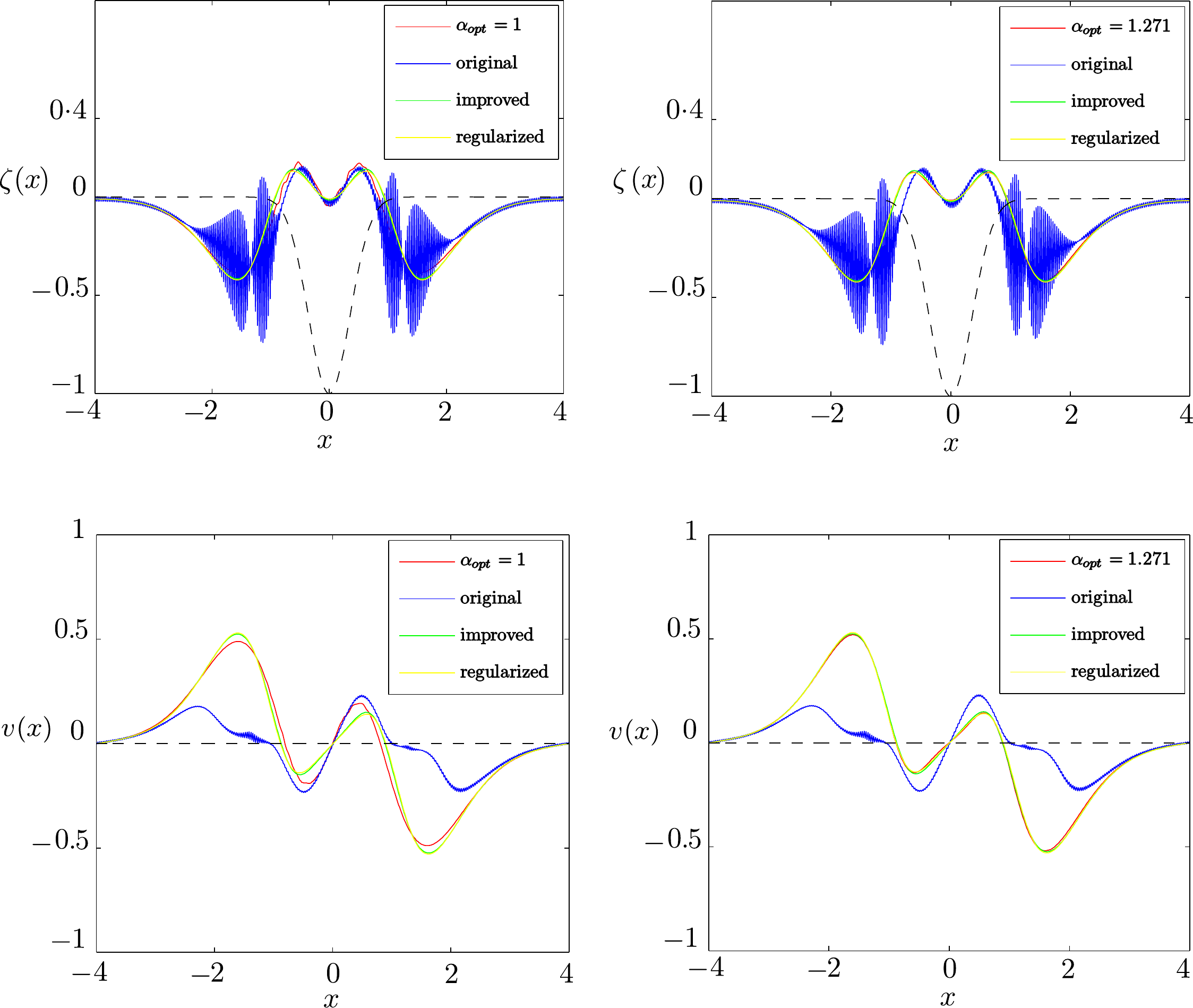}}
    \caption{Comparison with the Green-Naghdi models, without surface tension ($\bo^{-1}=0$), at time $t=2$, for $\alpha_{opt}=1$ (left) and $\alpha_{opt}=1.271$ (right)}
    \label{DIT3}
\end{figure} 
        \begin{figure}[H]
    \centering
{\includegraphics[scale=0.68]{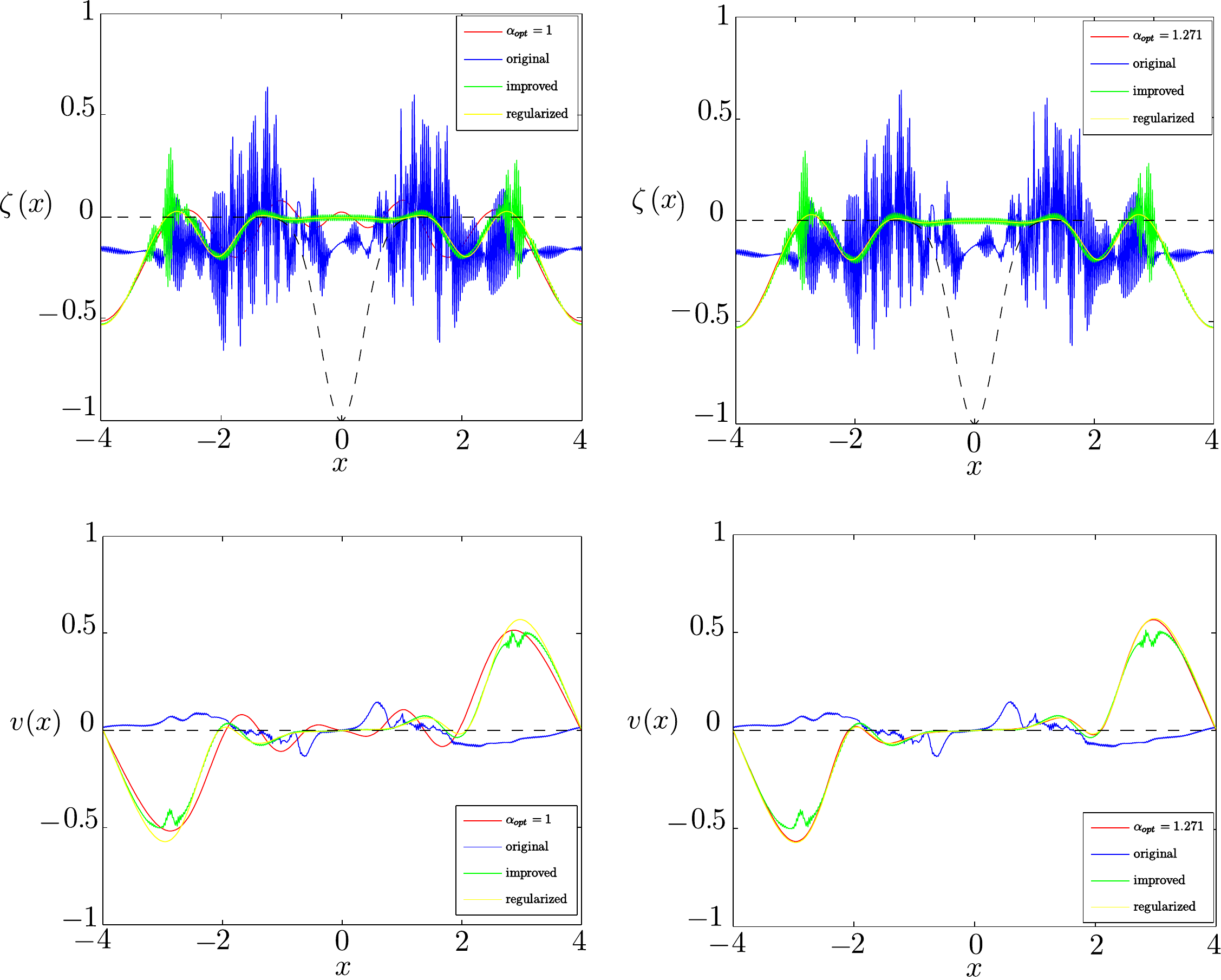}}
    \caption{Comparison with the Green-Naghdi models, without surface tension ($\bo^{-1}=0$), at time $t=5$, for $\alpha_{opt}=1$ (left) and $\alpha_{opt}=1.271$ (right)}
    \label{DIT4}
\end{figure} 
Figures~\ref{DIT1} and~\ref{DIT2} show the comparisons between our numerical solution for $\alpha_{opt}=1$ (left) and $\alpha_{opt}=1.271$ (right) and the Green-Naghdi models solutions obtained in~\cite{DucheneIsrawiTalhouk16}, with a small amount of surface tension, at time $t=2$ and $t=3$ respectively.
 We observe an excellent agreement between our numerical solution computed for $\alpha_{opt}=1.271$ and both ``improved" and ``regularized"  models at $t=2$ and $t=3$. As expected, at $t=3$ the original model induces Kelvin-Helmholtz instabilities. Meanwhile, the flows predicted by the regularized and improved models and by our model~\eqref{GNCH6} with $\alpha_{opt}=1.271$ remain smooth and are very similar. Similarly, Figure~\ref{DIT3} shows an excellent agreement between the numerical solutions computed for $\alpha_{opt}=1.271$ with the ``improved" and ``regularized"  models, without surface tension at time $t=2$, while the flow of the original model is completely destroyed due to Kelvin-Helmholtz instabilities. However, at a larger time ($t=5$), figure~\ref{DIT4} shows that in absence of surface tension ($\bo^{-1}=0$) Kelvin-Helmholtz instabilities will destroy completely both ``original" and ``improved"  models while the numerical solution computed for $\alpha_{opt}=1.271$ and for the ``regularized" models remain smooth and very similar. 
 
  The overall observations show the importance of the choice of the parameter $\alpha$ as in Section~\ref{Secalphachoice} in improving the frequency dispersion. Indeed, when choosing $\alpha_{opt}=1.271$, we observe an excellent matching between our numerical solutions and those obtained by the ``improved" model before the latter is completely destroyed in absence of surface tension due to the Kelvin-Helmholtz instabilities. As well, our numerical solution matches the one computed by the ``regularized" model even for a large time and with or without surface tension. This is not the case when choosing $\alpha_{opt}=1$. In fact, the ``improved" model has the same dispersion relation as the one of the \emph {full Euler system} and the dispersion relation of the ``regularized" model fit the one of the \emph {full Euler system} to an $\OO(\mu^3)$ order. This explains the reason behind the matching when choosing an optimal value for $\alpha$.
    \subsubsection{Dam-break problem in the two layers case}\label{DAM2}
\begin{figure}[H]
    \centering
{\includegraphics[scale=0.52]{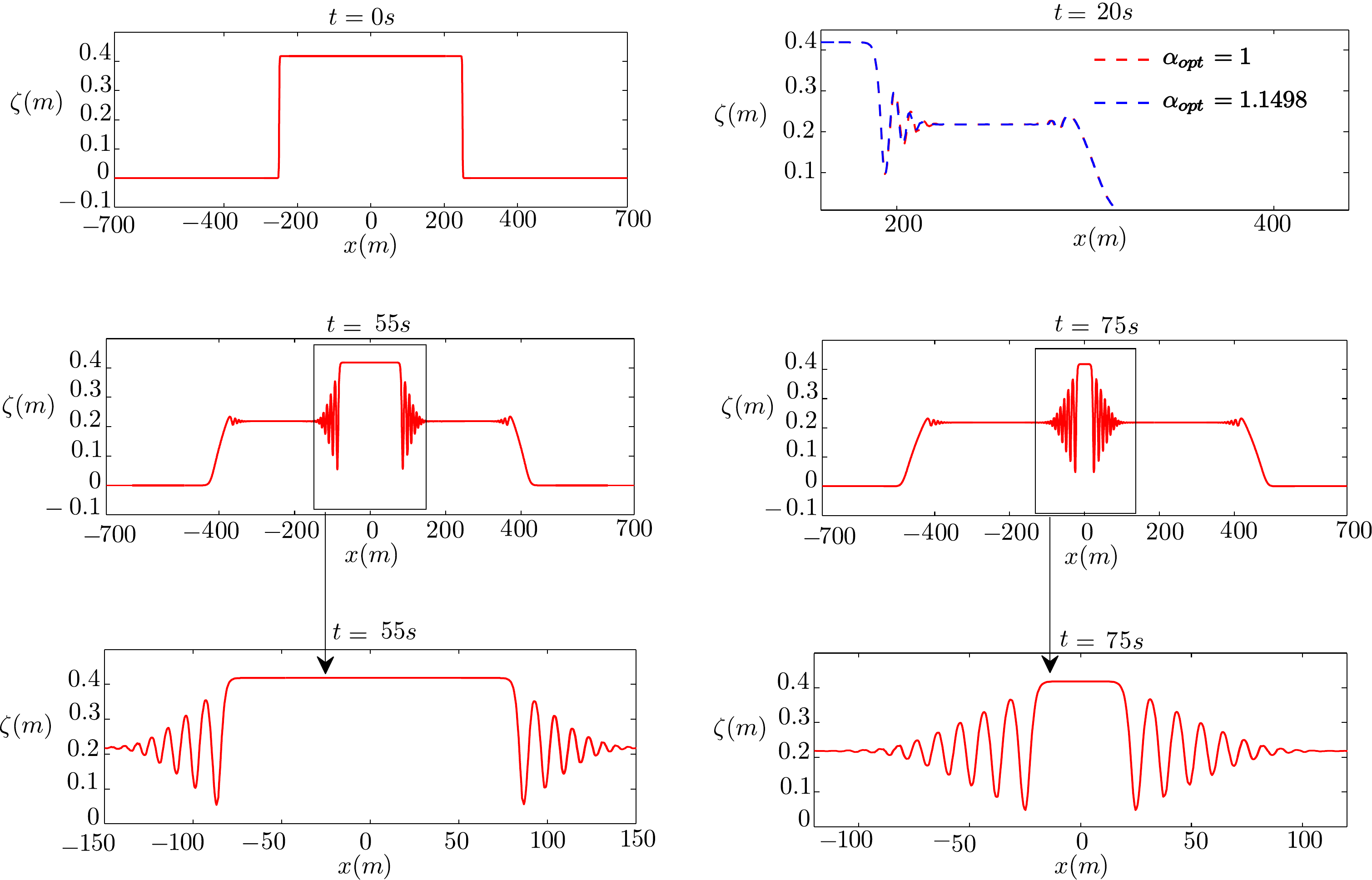}}
    \caption{Dam break in the two layers case: water surface profiles at several times}
    \label{dam2}
\end{figure}

This simulation concerns a test of the ability of our numerical scheme~\eqref{GNCH6dim} to handle discontinuities when considering a dam-break problem in the two layers case. To this end, we use the same initial data   as in the one layer case given by~\eqref{daminitialdata}, where $a=0.2091 \ m$. The simulations are performed on the interval $x \in (-700,700)$, discretized with 2800 cells imposing periodic conditions on each boundary. The dimensionless parameters representing the ratio between the depth and the ratio between the densities of the two layers are set respectively as follows: $\delta=0.5$, $\gamma=0.95$. Taking into account a small amount of surface tension, we set $\bo^{-1}= 5\times 10^{-5}$. We would like to mention that the simulations are computed using $\alpha=1$. As expected, the same simulations performed when choosing $\alpha=1.1498$ yielded the same result since $\alpha_{opt}\rightarrow 1$ when considering large wave numbers as explained in Section~\ref{Secalphachoice}. Figure~\ref{dam2} shows the results of the numerical simulation at several times, generating two dispersive tails propagating towards the center and two rarefaction waves that travel in opposite directions, on each side of the “dam”. A zoom on the dispersive tails is proposed at $t=55 \ s$ and $t=75 \ s$. Indeed, capturing this phenomenon accurately exhibit the high accuracy of our numerical scheme.
\section{Conclusion}
In this work, we presented a numerical scheme for the Green-Naghdi model in the Camassa-Holm regime. This model is first reformulated under more appropriate structure, where the time dependency of a second order differential operator present in the model is removed, keeping the stabilizing effects of its inversion. Furthermore, we improved the frequency dispersion of the original model keeping the same order of precision thanks to a parameter $\alpha$ to be precisely chosen. Additionally, our model do not contain third order derivatives that may create high frequencies instabilities. We then propose an efficient, precise and stable numerical splitting scheme that decomposes the hyperbolic and dispersive parts of the equations. The hyperbolic part is treated with a finite-volume method, where we consider the following schemes with different orders of accuracy: the first order finite-volume scheme based on a VFRoe method, the second order finite volume scheme following the classical MUSCL approach and finally a fifth-order WENO reconstruction. On the other hand, we treat the dispersive part with a finite-difference scheme, using second and fourth order formulas.  
Concerning time discretization we use classical second and fourth order Runge-Kutta methods, according to the order of the space derivative in consideration. Finally, we present several numerical validations in the one and two layers cases, showing the numerical efficiency and accuracy of our scheme and exhibiting its ability to reduce numerical dispersion and dissipation and to deal with discontinuities. The next step of this study may concern the extension of this numerical scheme to a more general configuration of variable topography. We believe that this splitting strategy may be applied in the variable bottom case.

%
%
\end{document}